\documentclass[a4paper,11pt]{article}

\usepackage[all]{xy}
\usepackage[T1]{fontenc}
\usepackage[dvips]{graphicx}
\usepackage{amsmath,amssymb,amsfonts,amsthm,latexsym,mathrsfs,textcomp,verbatim, enumerate, soul, color}
\xyoption{web}

\title{The unification of Mathematics\\ via\\ Topos Theory}
\author{Olivia Caramello\thanks{Supported by a visiting position from the Centro di Ricerca Matematica Ennio De Giorgi, Pisa.}}

\date{June 20, 2010}
\begin{document}

\mathcode`\<="4268  
\mathcode`\>="5269  
\mathcode`\.="313A  
\mathchardef\semicolon="603B 
\mathchardef\gt="313E
\mathchardef\lt="313C

\newcommand{\app}
 {{\sf app}}

\newcommand{\Ass}
 {{\bf Ass}}

\newcommand{\ASS}
 {{\mathbb A}{\sf ss}}

\newcommand{\Bb}
{\mathbb}

\newcommand{\biimp}
 {\!\Leftrightarrow\!}

\newcommand{\bim}
 {\rightarrowtail\kern-1em\twoheadrightarrow}

\newcommand{\bjg}
 {\mathrel{{\dashv}\,{\vdash}}}

\newcommand{\bstp}[3]
 {\mbox{$#1\! : #2 \bim #3$}}

\newcommand{\cat}
 {\!\mbox{\t{\ }}}

\newcommand{\cinf}
 {C^{\infty}}

\newcommand{\cinfrg}
 {\cinf\hy{\bf Rng}}

\newcommand{\cocomma}[2]
 {\mbox{$(#1\!\uparrow\!#2)$}}

\newcommand{\cod}
 {{\rm cod}}

\newcommand{\comma}[2]
 {\mbox{$(#1\!\downarrow\!#2)$}}

\newcommand{\comp}
 {\circ}

\newcommand{\cons}
 {{\sf cons}}

\newcommand{\Cont}
 {{\bf Cont}}

\newcommand{\ContE}
 {{\bf Cont}_{\cal E}}

\newcommand{\ContS}
 {{\bf Cont}_{\cal S}}

\newcommand{\cover}
 {-\!\!\triangleright\,}

\newcommand{\cstp}[3]
 {\mbox{$#1\! : #2 \cover #3$}}

\newcommand{\Dec}
 {{\rm Dec}}

\newcommand{\DEC}
 {{\mathbb D}{\sf ec}}

\newcommand{\den}[1]
 {[\![#1]\!]}

\newcommand{\Desc}
 {{\bf Desc}}

\newcommand{\dom}
 {{\rm dom}}

\newcommand{\Eff}
 {{\bf Eff}}

\newcommand{\EFF}
 {{\mathbb E}{\sf ff}}

\newcommand{\empstg}
 {[\,]}

\newcommand{\epi}
 {\twoheadrightarrow}

\newcommand{\estp}[3]
 {\mbox{$#1 \! : #2 \epi #3$}}

\newcommand{\ev}
 {{\rm ev}}

\newcommand{\Ext}
 {{\rm Ext}}

\newcommand{\fr}
 {\sf}

\newcommand{\fst}
 {{\sf fst}}

\newcommand{\fun}[2]
 {\mbox{$[#1\!\to\!#2]$}}

\newcommand{\funs}[2]
 {[#1\!\to\!#2]}

\newcommand{\Gl}
 {{\bf Gl}}

\newcommand{\hash}
 {\,\#\,}

\newcommand{\hy}
 {\mbox{-}}

\newcommand{\im}
 {{\rm im}}

\newcommand{\imp}
 {\!\Rightarrow\!}

\newcommand{\Ind}[1]
 {{\rm Ind}\hy #1}

\newcommand{\iten}[1]
{\item[{\rm (#1)}]}

\newcommand{\iter}
 {{\sf iter}}

\newcommand{\Kalg}
 {K\hy{\bf Alg}}

\newcommand{\llim}
 {{\mbox{$\lower.95ex\hbox{{\rm lim}}$}\atop{\scriptstyle
{\leftarrow}}}{}}

\newcommand{\llimd}
 {\lower0.37ex\hbox{$\pile{\lim \\ {\scriptstyle
\leftarrow}}$}{}}

\newcommand{\Mf}
 {{\bf Mf}}

\newcommand{\Mod}
 {{\bf Mod}}

\newcommand{\MOD}
{{\mathbb M}{\sf od}}

\newcommand{\mono}
 {\rightarrowtail}

\newcommand{\mor}
 {{\rm mor}}

\newcommand{\mstp}[3]
 {\mbox{$#1\! : #2 \mono #3$}}

\newcommand{\Mu}
 {{\rm M}}

\newcommand{\name}[1]
 {\mbox{$\ulcorner #1 \urcorner$}}

\newcommand{\names}[1]
 {\mbox{$\ulcorner$} #1 \mbox{$\urcorner$}}

\newcommand{\nml}
 {\triangleleft}

\newcommand{\ob}
 {{\rm ob}}

\newcommand{\op}
 {^{\rm op}}
 
\newcommand{\palrr}[4]{ 
  \def\labelstyle{\scriptstyle} 
  \xymatrix{ {#1} \ar@<0.5ex>[r]^{#2} \ar@<-0.5ex>[r]_{#3} & {#4} } } 
  
\newcommand{\palrl}[4]{ 
  \def\labelstyle{\scriptstyle} 
  \xymatrix{ {#1} \ar@<0.5ex>[r]^{#2}  &  \ar@<0.5ex>[l]^{#3} {#4} } }  

\newcommand{\pepi}
 {\rightharpoondown\kern-0.9em\rightharpoondown}

\newcommand{\pmap}
 {\rightharpoondown}

\newcommand{\Pos}
 {{\bf Pos}}

\newcommand{\prarr}
 {\rightrightarrows}

\newcommand{\princfil}[1]
 {\mbox{$\uparrow\!(#1)$}}

\newcommand{\princid}[1]
 {\mbox{$\downarrow\!(#1)$}}

\newcommand{\prstp}[3]
 {\mbox{$#1\! : #2 \prarr #3$}}

\newcommand{\pstp}[3]
 {\mbox{$#1\! : #2 \pmap #3$}}

\newcommand{\relarr}
 {\looparrowright}

\newcommand{\rlim}
 {{\mbox{$\lower.95ex\hbox{{\rm lim}}$}\atop{\scriptstyle
{\rightarrow}}}{}}

\newcommand{\rlimd}
 {\lower0.37ex\hbox{$\pile{\lim \\ {\scriptstyle
\rightarrow}}$}{}}

\newcommand{\rstp}[3]
 {\mbox{$#1\! : #2 \relarr #3$}}

\newcommand{\scn}
 {{\bf scn}}

\newcommand{\scnS}
 {{\bf scn}_{\cal S}}

\newcommand{\semid}
 {\rtimes}

\newcommand{\Sep}
 {{\bf Sep}}

\newcommand{\sep}
 {{\bf sep}}

\newcommand{\Set}
 {{\bf Set}}

\newcommand{\Sh}
 {{\bf Sh}}

\newcommand{\ShE}
 {{\bf Sh}_{\cal E}}

\newcommand{\ShS}
 {{\bf Sh}_{\cal S}}

\newcommand{\sh}
 {{\bf sh}}

\newcommand{\Simp}
 {{\bf \Delta}}

\newcommand{\snd}
 {{\sf snd}}

\newcommand{\stg}[1]
 {\vec{#1}}

\newcommand{\stp}[3]
 {\mbox{$#1\! : #2 \to #3$}}

\newcommand{\Sub}
 {{\rm Sub}}

\newcommand{\SUB}
 {{\mathbb S}{\sf ub}}

\newcommand{\tbel}
 {\prec\!\prec}

\newcommand{\tic}[2]
 {\mbox{$#1\!.\!#2$}}

\newcommand{\tp}
 {\!:}

\newcommand{\tps}
 {:}

\newcommand{\tsub}
 {\pile{\lower0.5ex\hbox{.} \\ -}}

\newcommand{\wavy}
 {\leadsto}

\newcommand{\wavydown}
 {\,{\mbox{\raise.2ex\hbox{\hbox{$\wr$}
\kern-.73em{\lower.5ex\hbox{$\scriptstyle{\vee}$}}}}}\,}

\newcommand{\wbel}
 {\lt\!\lt}

\newcommand{\wstp}[3]
 {\mbox{$#1\!: #2 \wavy #3$}}
 
\newcommand{\fu}[2]
{[#1,#2]}


%
%
%
\def\pushright#1{{
   \parfillskip=0pt            
   \widowpenalty=10000         
   \displaywidowpenalty=10000  
   \finalhyphendemerits=0      
  %
   \leavevmode                 
   \unskip                     
   \nobreak                    
   \hfil                       
   \penalty50                  
   \hskip.2em                  
   \null                       
   \hfill                      
   {#1}                        
  %
   \par}}                      

\def\qed{\pushright{$\square$}\penalty-700 \smallskip}

\newtheorem{theorem}{Theorem}[section]

\newtheorem{proposition}[theorem]{Proposition}

\newtheorem{scholium}[theorem]{Scholium}

\newtheorem{lemma}[theorem]{Lemma}

\newtheorem{corollary}[theorem]{Corollary}

\newtheorem{conjecture}[theorem]{Conjecture}

\newenvironment{proofs}%
 {\begin{trivlist}\item[]{\bf Proof }}%
 {\qed\end{trivlist}}

  \newtheorem{rmk}[theorem]{Remark}
\newenvironment{remark}{\begin{rmk}\em}{\end{rmk}}

  \newtheorem{rmks}[theorem]{Remarks}
\newenvironment{remarks}{\begin{rmks}\em}{\end{rmks}}

  \newtheorem{defn}[theorem]{Definition}
\newenvironment{definition}{\begin{defn}\em}{\end{defn}}

  \newtheorem{eg}[theorem]{Example}
\newenvironment{example}{\begin{eg}\em}{\end{eg}}

  \newtheorem{egs}[theorem]{Examples}
\newenvironment{examples}{\begin{egs}\em}{\end{egs}}


\bgroup           
\let\footnoterule\relax  
\maketitle

\begin{abstract}
We present a set of principles and methodologies which may serve as foundations of a unifying theory of Mathematics. These principles are based on a new view of Grothendieck toposes as unifying spaces being able to act as `bridges' for transferring information, ideas and results between distinct mathematical theories. 
\end{abstract} 
\egroup

\vspace{6.5cm}
{\small \textsc{Centro di Ricerca Matematica Ennio De Giorgi, Scuola Normale Superiore, Piazza dei Cavalieri 3, 56100 Pisa, Italy}\\
\emph{E-mail address:} \texttt{olivia.caramello@sns.it}}

\vspace{0.5cm}

{\small \textsc{Department of Pure Mathematics and Mathematical Statistics, University of Cambridge, Wilberforce Road, Cambridge CB3 0WB, UK}\\
\emph{E-mail address:} \texttt{O.Caramello@dpmms.cam.ac.uk}}

\vspace{5 mm}

\tableofcontents

\newpage

\section{Introduction}

In this paper, our aim is to lay a set of principles and methodologies which may serve as foundations of a unifying theory of Mathematics. By this, we mean a meta-theory offering methodologies for comparing different mathematical theories with each other, discovering analogies between them as well as pointing out peculiarities, and, most importantly, providing effective means for transferring results and techniques between distinct fields. 

Mathematics is divided into several distinct areas: Algebra, Analysis, Geometry, Topology, Number Theory etc. Each of these areas has evolved throughout the years by developing its own ideas and techniques, and has reached by now a remarkable degree of specialization. With time, various connections between the areas have been discovered, leading in some cases to the creation of actual `bridges' between different mathematical branches; many times methods of one field have been employed to derive results in another, and this interplay of different points of view in a same subject has always had a fundamental role in illuminating the nature of concepts, establishing new results and suggesting new lines of investigation.
 
The profound reasons behind the analogies and connections between different fields are best understood in the context of Mathematical Logic; indeed, Logic provides a means for formalizing any kind of mathematical concept, so that the investigation of the relationships between different theories can be carried out at a completely rigorous level.

Considered the importance of building bridges between distinct mathematical branches, it would be highly desirable that Logic would not just serve as a tool for analyzing analogies already discovered in Mathematics but could instead play an active role in identifying new connections between existing fields, as well as suggesting new directions of mathematical investigation. As it happens, we now have enough mathematical tools at our disposal for trying to achieve this goal.

By providing a system in which all the usual mathematical concepts can be expressed rigorously, Set Theory has represented the first serious attempt of Logic to unify Mathematics at least at the level of language. Later, Category Theory offered an alternative abstract language in which most of Mathematics can be formulated and, as such, has represented a further advancement towards the goal of `unifying Mathematics'. Anyway, both these systems realize a unification which is still limited in scope, in the sense that, even though each of them provides a way of expressing and organizing Mathematics in one single language, they do not offer by themselves effective methods for an actual transfer of knowledge between distinct fields. On the other hand, the principles that we will sketch in the present paper define a different and more substantial approach to the unification of Mathematics. 

Our methodologies are based on a new view of Grothendieck toposes as unifying spaces which can serve as `bridges' for transferring information, ideas and results between distinct mathematical theories. This view of toposes has emerged in the course of the Ph.D. investigations carried out by the author at the University of Cambridge from 2006 to 2009, and in fact the results in \cite{OC} provide compelling technical evidence for the validity of such view (indeed, it is also among the purposes of the present paper to serve as a `conceptual guide' to the methods in the author's dissertation).

In this paper, I give an outline of the fundamental principles characterizing my view of toposes as `bridges' connecting different mathematical theories and describe the general methodologies which have arisen from such a view and which have motivated my investigations so far. The analysis will be complemented by the discussion of the results in \cite{OC} which best illustrate the application of the above-mentioned principles. These principles are abstract and transversal to the various mathematical fields, and the application of them can lead to a huge amount of surprising and non-trivial results in any area of Mathematics, so we hope that the reader will get motivated to try out these methods in his or her fields of interest.

\section{Geometric theories and their classifying toposes}

In this section we discuss the background on geometric theories and classifying toposes which is necessary for the following parts of the paper (for a comprehensive source of information on the subject we refer the reader to Part D of \cite{El}). 

\subsection{Geometric theories}

\emph{Geometric theories} are a wide class of (multi-sorted) infinitary first-order theories. Recall that a theory over a signature $\Sigma$ is geometric if its axioms can be presented in the form $\forall \vec{x}(\phi \imp \psi)$, where $\phi$ and $\psi$ are geometric formulae over $\Sigma$ (i.e. formulae with a finite number of free variables built up from atomic formulae over $\Sigma$ by only using finitary conjunctions, infinitary disjunctions, and existential quantifications) in the context $\vec{x}$. The axioms of geometric theories are often presented in sequent form, that is one writes $\phi \vdash_{\vec{x}} \psi$ instead of $(\forall \vec{x})(\phi \imp \psi)$. If in the above-mentioned definition of geometric theory we replace `infinitary disjunctions' with `finitary disjunctions', we get the class of theories known as \emph{coherent theories}.

The attribute `geometric' should not induce the reader to think that this is a class of theories which have something particular to do with Geometry (apart from the first historical applications); the concept is completely general and geometric theories can be found in essentially every field of Mathematics.  
 
Indeed, it is a matter of fact that most of the theories naturally arising in Mathematics have a geometric axiomatization (over their signature). Anyway, if a finitary first-order theory $\mathbb T$ is not geometric, we can canonically construct a coherent theory over a larger signature, called the \emph{Morleyization} of $\mathbb T$, whose models in the category $\Set$ of sets (more generally, in any Boolean coherent category) can be identified with those of $\mathbb T$ (cfr. Lemma D1.5.13 \cite{El}).

The notion of Morleyization is important because it enables us to study any kind of first-order theory by using the methods of Topos Theory. In fact, we can expect many important properties of first-order theories to be naturally expressible as properties of their Morleyizations, and these latter properties to be in turn expressible in terms of `invariants' of their classifying toposes. For example, a first-order theory is complete if and only if its Morleyization ${\mathbb T}'$ satisfies the following property: any coherent (equivalently, geometric) sentence over its signature is provably equivalent to $\bot$ or to $\top$, but not both; and this property is precisely equivalent to saying that its classifying topos is two-valued. 

Due to its infinitary nature, geometric logic is quite expressive; important mathematical properties which are not expressible in finitary first-order logic (over a given signature) often admit a geometric axiomatization. For example, the property of an element of a commutative ring with unit to be nilpotent is well-known not to be expressible by a finitary first-order formula over the signature of rings, but it is obviously expressible as a geometric formula over this signature. An interesting infinitary geometric theory, which was studied topos-theoretically in \cite{OC4}, is the theory of fields of finite characteristic which are algebraic over their prime fields. Another example of a geometric theory which is not finitary first-order is given by the theory of torsion abelian groups. 

Geometric theories, including the infinitary ones, should definitely be regarded as objects which are worth of investigation. The infinitary nature of geometric logic has surely represented an obstacle towards the recognition of this class of theories as a subject of primary importance; in fact, one cannot in general employ the methods of classical model theory for studying geometric theories, and therefore methods of completely different nature must be used.   

As we shall see in section \ref{logGrothendieck}, from the point of view of Topos Theory the finiteness assumptions on the nature of the logic indeed appear to be very unnatural restrictions; even starting with finitary theories, studying them topos-theoretically often leads to the consideration of infinitary ones. The logic underlying Grothendieck toposes is geometric logic, in its full infinitary nature (cfr. section \ref{logGrothendieck}); and, as we shall argue below, Topos Theory provides a set of incredibly powerful techniques for studying geometric theories. 

We emphasize that geometric theories, as any kind of first-order theories, are objects of purely syntactic nature. As in classical finitary first-order logic we have a syntactic notion of provability of first-order sentences (relative to a theory), so we have a natural proof system for geometric (resp. coherent) logic, described in terms of inference rules involving geometric (resp. coherent) sequents, which yields a notion of provability of geometric (resp. coherent) sequents relative to a given geometric (resp. coherent) theory. In geometric logic, classical and intuitionistic provability of geometric sequents coincide, so we can well omit the law of excluded middle from these proof systems without affecting the corresponding notion of provability. We refer the reader to Part D of \cite{El} for a detailed presentation of these systems. 

It is well-known that first-order languages can always be interpreted in the context of (a given model of) set theory. In fact, these languages can also be meaningfully interpreted in a category, provided that the latter possesses enough categorical structure to allow the interpretation of the connectives and quantifiers arising in the formulae of the language (in this categorical semantics sorts are interpreted as objects, terms as arrows and formulae as subobjects, in a way that respects the logical structure of compound expressions); this leads, under the same assumptions, to a notion of satisfaction of sequents in a categorical structure, and hence to a notion of model of a first-order theory in a category, which specializes, in the case of the category of sets, to the classical Tarskian definition of (set-based) model of a first-order theory. For example, a topological group can be seen as a model of the theory of groups in the category of topological spaces, and a sheaf of rings on a topological space $X$ can be seen as a model of the theory of rings in the topos $\Sh(X)$. This categorical semantics, when it is defined, is always sound with respect to the above-mentioned proof systems, provided that the latter do not contain the law of excluded middle. 

Grothendieck toposes are a wide class of categories in which we can interpret geometric logic, and in fact they provide us with a strong form of completeness for geometric theories. Indeed, for any geometric theory $\mathbb T$, there is a Grothendieck topos $\Set[{\mathbb T}]$, namely the \emph{classifying topos} of $\mathbb T$, which contains a \emph{conservative} model of $\mathbb T$ (i.e. a model $U_{\mathbb T}$ such that the geometric sequents provable in $\mathbb T$ are exactly those which are satisfied in $U_{\mathbb T}$), which is moreover \emph{universal} in the sense that any model of $\mathbb T$ in a Grothendieck topos arises, up to isomorphism, as the image of $U_{\mathbb T}$ under the inverse image functor of a (unique up to equivalence) geometric morphism of toposes from the topos in which the model lives to the classifying topos (cfr. section \ref{classif} below). We remark that, due to their infinitary nature, geometric theories do not in general enjoy a classical form of completeness i.e. it is not in general true that a sequent which is valid in all the set-based models of a geometric theory is provable in the theory by using geometric logic. Anyway, as we have said, the concept of universal model of a geometric theory yields a strong form of completeness, and in fact we believe that it should be taken as the fundamental tool (replacing the standard ones) in the investigation of the aspects of completeness of the theory (cfr. sections \ref{classif} and \ref{logGrothendieck}). 

\subsection{Classifying toposes}\label{classif}

A \emph{classifying topos} of a geometric theory $\mathbb T$ over a signature $\Sigma$ is a\\ Grothendieck topos $\Set[{\mathbb T}]$ such that for any Grothendieck topos $\cal E$ the category ${\bf Geom}({\cal E}, \Set[{\mathbb T}])$ of geometric morphisms from ${\cal E}$ to $\Set[{\mathbb T}]$ is equivalent to the category of models of $\mathbb T$ in the topos $\cal E$, naturally in $\cal E$; naturality means that for any geometric morphism $f:{\cal E}\to {\cal F}$ of Grothendieck toposes, we have a commutative square
\[  
\xymatrix {
{\bf Geom}({\cal F}, \Set[{\mathbb T}]) \ar[d]^{- \circ f} \ar[rr]^{\simeq}   & & {{\mathbb T}}{\textrm{-mod}}({\cal F})  \ar[d]^{{\mathbb T}\textbf{-mod}(f^{\ast})}  \\
{\bf Geom}({\cal E}, \Set[{\mathbb T}])  \ar[rr]^{\simeq} & &   {\mathbb T}{\textrm{-mod}}({\cal E})}
\]  
in the (meta-)category of categories $\textbf{CAT}$.

Recall that geometric morphisms are the natural, topologically motivated, notion of morphism between Grothendieck toposes; indeed, the assignment sending a locale $L$ to the topos of sheaves $\Sh(L)$ on it gives rise to a full and faithful functor from the category of locales to the category of Grothendieck toposes and geometric morphisms between them, which identifies the former category as a full reflective subcategory of the latter. 

The classifying topos of a geometric theory $\mathbb T$ can be seen as a \emph{representing object} for the $\mathbb T$-model (pseudo-)functor
\[
{{\mathbb T}\textrm{-mod}:\mathfrak{BTop}^{\textrm{op}} \to \textbf{CAT}}
\]
from the opposite of the category $\mathfrak{BTop}$ of Grothendieck toposes to the (meta-)category of categories $\textbf{CAT}$ which assigns 
\begin{itemize}
\item to a topos $\cal E$ the category ${\mathbb T}\textbf{-mod}({\cal E})$ of models of $\mathbb T$ in $\cal E$ and 
\item to a geometric morphism $f:{\cal E}\to {\cal F}$ the functor\\ ${\mathbb T}\textbf{-mod}(f^{\ast}):{\mathbb T}\textbf{-mod}({\cal F}) \to {\mathbb T}\textbf{-mod}({\cal E})$\\ sending a model $M\in {\mathbb T}\textbf{-mod}({\cal F})$ to its image $f^{\ast}(M)$ under the inverse image functor $f^{\ast}$ of $f$.
\end{itemize} 
In particular, classifying toposes are \emph{unique up to categorical equivalence}.

A concept which is strictly connected with that of classifying topos is that of universal model. A \emph{universal model} of a geometric theory $\mathbb T$ is a model \emph{$U_{\mathbb T}$} of $\mathbb T$ in a Grothendieck topos $\cal G$ such that for any $\mathbb T$-model $M$ in a Grothendieck topos $\cal F$ there exists a unique (up to isomorphism) geometric morphism $f_{M}:{\cal F}\to {\cal G}$ such that $f_{M}^{\ast}(U_{\mathbb T})\cong M$.

By the ($2$-dimensional) Yoneda Lemma, if a topos $\cal G$ contains a \emph{universal model} of a geometric theory $\mathbb T$ then $\cal G$ satisfies the universal property of the \emph{classifying topos} of $\mathbb T$. Conversely, if a topos $\cal E$ classifies a geometric theory $\mathbb T$ then $\cal E$ contains a universal model of $\mathbb T$.

Therefore universal models, as well as classifying toposes, are uniquely determined up to categorical equivalence. Specifically, if $M$ and $N$ are universal models of a geometric theory $\mathbb T$ lying respectively in toposes $\cal F$ and $\cal G$ then there exists a unique (up to isomorphism) geometric equivalence between $\cal F$ and $\cal G$ such that its inverse image functors send $M$ and $N$ to each other (up to isomorphism).
 
Classifying toposes of arbitrary first-order theories do not in general exist, for the essential reason that inverse image functors of geometric morphisms may not preserve the interpretation of implications and universal quantifications (even though they always preserve the interpretation of geometric formulae). On the other hand, classifying toposes always exist for geometric theories, and they can be built canonically from them by means of a syntactic construction (cfr. section \ref{onetopos} below). This fact is of fundamental importance, at least for our purposes, and we shall extensively discuss it below. 

We remark that the mere \emph{existence} of classifying toposes is alone a fact of primary conceptual importance; the fact that the $\mathbb T$-model functor is representable means that, in a sense, all the information about the functor or, in other words, about the categories of models of the theory $\mathbb T$, is condensed in one single object, namely the classifying topos. This clearly represents a `symmetry' result within the model theory of geometric theories, which demonstrates that the environment of Grothendieck toposes enjoys a strong form of `inner completeness' with respect to geometric theories: all the models of a given geometric theory in Grothendieck toposes, including all the classical set-based models of the theory, are images of a single `universal model' lying in its classifying topos, in which the integration of syntactic and semantic aspects of theory takes place. If we restrict our considerations to the set-based models of a theory $\mathbb T$, we do not have a similar representability result: we have to expand our view, for example to the world of Grothendieck toposes, in order to find this kind of `symmetry'. After all, this kind of phenomenon is quite common in Mathematics; think for example that, before the invention of the complex plane, analysts were having a hard time trying to understand the behaviour of solutions to polynomial equations in the real line, while as soon as `the natural environment' for polynomial equations was eventually found, everything started to be perceived under a new light, and people immediately stopped wondering why the behaviour of polynomial equations in the real line appeared so erratic. In fact, the fundamental theorem of algebra perfectly embodies this achieved `symmetry' with respect to polynomial equations.  

A familiar image that comes to mind in thinking about this is that of a \emph{sun} and the \emph{shadows} that is generates. Shadows arise when the light of the sun meets some form of solid matter; similarly, models arise when a piece of syntax is interpreted in a given `concrete' environment. For example, the abstract, syntactic, notion of group gives rise to many different models in different categories, namely the classical notion of group (if we interpret it in the category of sets), the notion of topological group (if we interpret it in the category of topological spaces), the notion of algebraic group (if we interpret it in the category of algebraic varieties), the notion of Lie group (if we interpret it in the category of smooth manifolds), etc.

As shadows are easily comparable with each other, if one bears in mind that they ultimately come from a unique source, so the study of models of a theory can be very much aided by the consideration of the syntactic aspects of the theory. For example, all the abstract algebraic properties of groups (in the sense of algebraic identities which are provable in the axiomatic theory of groups) can be canonically interpreted in all of the above-mentioned categories to `automatically' yield new results on topological groups, algebraic groups, Lie groups etc. All of these notions share a common core which lies at the syntactic level rather than at the semantical one.

As the sun is a unifying source of shadows, syntax plays the role of a unifying concept among mathematical structures. This is perhaps not particularly visible in the traditional setting of classical (first-order) model theory, where we have a completeness theorem for finitary first-order logic providing an equivalence between the syntactic viewpoint and the semantical one based on set theory. But when one has many semantics available, as in topos-based model theory, syntax clearly acquires a prominent role. 

\section{The logic underlying Grothendieck toposes}\label{logGrothendieck}

We know that every geometric theory over a given signature has a classifying topos, which is uniquely determined up to categorical equivalence. This naturally leads to wondering whether this assignment is a sort of injection or surjection, that is to raising the following questions:\\ 

(1) Can two different geometric theories have equivalent classifying\\ toposes?\\

(2) Is every Grothendieck topos the classifying topos of some geometric theory?\\

The answers to both these questions are well-known. Specifically, it immediately follows from the characterization of the classifying topos as a representing object for the pseudofunctor of models (cfr. section \ref{classif}), that two geometric theories have equivalent classifying toposes if and only if they have equivalent categories of models in every Grothendieck topos $\cal E$, naturally in $\cal E$. Two such theories are said to be \emph{Morita-equivalent}. (We will come back to the subject of Morita-equivalence in section \ref{morita} below.) This answers the first question.

To answer the second question, we need to go slightly deeper into the structure of Grothendieck toposes. By definition, a Grothendieck topos is any (category equivalent to a) category $\Sh({\cal C}, J)$ of sheaves on a site $({\cal C}, J)$ (a site of definition of a Grothendieck topos $\cal E$ is any site $({\cal C}, J)$ such that the category $\Sh({\cal C}, J)$ of sheaves on $({\cal C}, J)$ is equivalent to $\cal E$); recall that a site $({\cal C}, J)$ consists of a small category $\cal C$ and a Grothendieck topology $J$ on $\cal C$ (we refer the reader to \cite{MM} for a first introduction to Topos Theory). 

Notice that there is an element of `non-canonicity' here, in that we cannot expect the topos $\Sh({\cal C}, J)$ to uniquely determine the site $({\cal C}, J)$; however, this aspect of `non-canonicity' is in many respects not at all an undesirable feature, and in fact it represents a fundamental ingredient of the view of toposes as unifying spaces described below.)

Recall that, given a site $({\cal C}, J)$, for every Grothendieck topos $\cal E$ we have an equivalence between the category ${\bf Geom}({\cal E}, \Sh({\cal C}, J))$ of geometric morphisms from $\cal E$ to $\Sh({\cal C}, J))$ and the category ${\bf Flat}_{J}(\cal{C}, \cal{E})$ of $J$-continuous flat functors from ${\cal{C}}$ to ${\cal{E}}$, naturally in $\cal E$. Now, we can construct a geometric theory ${\mathbb T}^{\cal C}_{J}$ such that its models in any Grothendieck topos $\cal E$ can be identified precisely with the $J$-continuous flat functors from ${\cal{C}}$ to ${\cal{E}}$ (and the homomorphisms of ${\mathbb T}^{\cal C}_{J}$-models can be identified with natural transformations between the corresponding flat functors); clearly, ${\mathbb T}^{\cal C}_{J}$ will be classified by the topos $\Sh({\cal C}, J)$. We call such theory ${\mathbb T}^{\cal C}_{J}$ the \emph{theory of $J$-continuous flat functors on $\cal C$}. This ensures that every Grothendieck topos arises as the classifying topos of some geometric theory, answering positively to the second question. It is instructive to write down explicitly an axiomatization of the theory ${\mathbb T}^{\cal C}_{J}$.  

The signature of ${\mathbb T}^{\cal C}_{J}$ has one sort $\name{A}$ for each object $A$ of $\cal C$, and one function symbol $\name{f}:\name{A} \to \name{B}$ for each arrow $f:A \to B$ in $\cal C$. The axioms of ${\mathbb T}^{\cal C}_{J}$ are the following (to indicate that a variable $x$ has sort $\name{A}$ we write $x^{A}$): 
\begin{equation}
(\top \vdash_{x} (\name{f}(x)=x))
\end{equation}

for any identity arrow $f$ in $\cal C$;
\begin{equation}
(\top \vdash_{x} (\name{f}(x)=\name{h}(\name{g}(x))))
\end{equation}
for any triple of arrows $f,g,h$ of $\cal C$ such that $f$ is equal to the composite $h\circ g$;
\begin{equation}
\top \vdash_{[]} \mathbin{\mathop{\textrm{\huge $\vee$}}\limits_{A\in Ob({\cal C})}}(\exists x^{A})\top
\end{equation}
(where the disjunction ranges over all the objects of $\cal C$);
\begin{equation}
(\top \vdash_{x^{A}, y^{B}} \mathbin{\mathop{\textrm{\huge $\vee$}}\limits_{A\stackrel{f}{\leftarrow} C \stackrel{g}{\rightarrow} B}}(\exists z^{C})(\name{f}(z^{C})=x^{A} \wedge \name{g}(z^{C})=y^{B}))
\end{equation}
for any objects $A$, $B$ of $\cal C$ (where the disjunction ranges over all the cones $A\stackrel{f}{\leftarrow} C \stackrel{g}{\rightarrow} B$ on the discrete diagram given by the pair of objects $A$ and $B$);
\begin{equation}
(\name{f}(x^{A})=\name{g}(x^{A}) \vdash_{x^{A}} \mathbin{\mathop{\textrm{\huge $\vee$}}\limits_{h:C\to A\in Eq(f,g)}}(\exists z^{C})(\name{h}(z^{C})=x^{A}))
\end{equation}
for any pair of arrows $f,g:A\to B$ in $\cal C$ with common domain and codomain (where the disjunction ranges over all the arrows $h$ which equalize $f$ and $g$);
\begin{equation}
(\top \vdash_{x^{A}} \mathbin{\mathop{\textrm{\huge $\vee$}}\limits_{i\in I}}(\exists y_{i}^{B_{i}})(\name{f_{i}}(y_{i}^{B_{i}})=x^{A}))  
\end{equation}
for each $J$-covering family $(f_{i}:B_{i}\to A \textrm{ | } i\in I)$.
          
Notice that the first two groups of axioms express functoriality, the third, fourth and fifth together express flatness (in terms of filteredness of the corresponding category of elements), while the sixth group of axioms expresses $J$-continuity. 

We remark that, even when $\cal C$ is cartesian, in which case flat functors on $\cal C$ have a finitary (coherent) axiomatization, the presence of the Grothendieck topology $J$ makes the axiomatization of ${\mathbb T}^{\cal C}_{J}$ in general infinitary. This shows that the `underlying logic' of Grothendieck toposes is indeed \emph{geometric logic}, in its full infinitary nature. Further evidence for this is provided by the duality theorem (cfr. section \ref{duality}), which asserts that the subtoposes of the classifying topos of a geometric theory correspond bijectively with the \emph{geometric} quotients of the theory. In fact, as we shall see below, most of the classical topos-theoretic invariants or constructions behave quite naturally with respect to geometric logic (that is, they correspond to natural properties or operations on geometric theories), but rather awkwardly with respect to the \emph{finitary} fragment of it i.e. coherent logic. Indeed, many important topos-theoretic constructions yield, even if applied to coherent theories, theories which are no longer coherent, a clear indication that these topos-theoretically induced transformations of theories can be naturally studied only in the full infinitary framework. Another illustration of the intrinsically geometric nature of the logic underlying Grothendieck toposes is given by the consideration of universal models of geometric theories; by Theorem 7.1.4 \cite{OC}, the subobjects (of an underlying object of) the universal model of a geometric theory can be identified with the (provable-equivalence classes of) \emph{geometric} formulae in a given context over the signature of the theory.     

In this connection, it seems illuminating also to recall the intrinsic characterization of geometric logic obtained in \cite{OC9}:

\begin{theorem}\label{teofond}

Let $\Sigma$ be a signature and $\cal S$ be a collection of $\Sigma$-structures in Grothendieck toposes closed under isomorphisms of structures. Then $\cal S$ is the collection of all models in Grothendieck toposes of a geometric theory over $\Sigma$ if and only if it satisfies the following two conditions:
\begin{enumerate}[(i)]
\item for any geometric morphism $f:{\cal F}\to {\cal E}$, if $M$ is in $\cal S$ then $f^{\ast}(M)$ is in $\cal S$;
\item for any (set-indexed) jointly surjective family $\{f_{i}:{\cal E}_{i}\to {\cal E} \textrm{ | } i\in I\}$ of geometric morphisms and any $\Sigma$-structure $M$ in $\cal E$, if $f_{i}^{\ast}(M)$ is in $\cal S$ for every $i\in I$ then $M$ is in $\cal S$. 
\end{enumerate}
\end{theorem}

\section{One topos, many sites}\label{onetopos}

The classifying topos of a geometric theory can always be built `canonically' from the theory by means of a syntactic construction: specifically, the classifying topos of a theory $\mathbb T$ is given by the category of sheaves on the \emph{geometric syntactic site} $({\cal C}_{\mathbb T}, J_{\mathbb T})$ of $\mathbb T$. Recall that the \emph{geometric syntactic category} ${\cal C}_{\mathbb T}$ of a geometric theory $\mathbb T$ over a signature $\Sigma$ has as objects the `renaming-equivalence classes' of geometric formulae-in-context $\{\vec{x}. \phi\}$ over $\Sigma$ and as arrows $[\theta]:\{\vec{x}. \phi\} \to \{\vec{y}. \psi\}$ the $\mathbb T$-provable equivalence classes of geometric formulae $[\theta]$ which are $\mathbb T$-provably functional from $\{\vec{x}. \phi\}$ to $\{\vec{y}. \psi\}$, while the \emph{geometric syntactic topology} $J_{\mathbb T}$ of $\mathbb T$ is the canonical Grothendieck topology on the geometric category ${\cal C}_{\mathbb T}$ (cfr. Part D of \cite{El} for more details).

If the theory $\mathbb T$ lies in a smaller fragment of geometric logic, such as cartesian, regular or coherent logic, the classifying topos can be alternatively calculated by taking the category of sheaves on other syntactic sites, namely the cartesian syntactic category ${\cal C}_{\mathbb T}^{\textrm{cart}}$ of $\mathbb T$ equipped with the trivial Grothendieck topology on it (if $\mathbb T$ is cartesian), the regular syntactic site $({\cal C}_{\mathbb T}^{\textrm{reg}}, J_{{\cal C}_{\mathbb T}^{\textrm{reg}}})$ (if $\mathbb T$ is regular), and the coherent syntactic site $({\cal C}_{\mathbb T}^{\textrm{coh}}, J_{{\cal C}_{\mathbb T}^{\textrm{coh}}})$ (if $\mathbb T$ is coherent).

So, for a given theory, there may be alternative ways of calculating its classifying topos, even remaining in the context of syntactic sites: for example, if a theory is cartesian, by regarding it as a cartesian, regular, coherent and geometric theory, we obtain four different `syntactic sites' such that the category of sheaves on them yields the classifying topos of the theory. This fact will be exploited in section \ref{examp} to derive various results in Logic. 

On the other hand, there is an alternative general method, of semantical nature, for calculating classifying toposes of geometric theories, based on the notion of theory of presheaf type. A theory is said to be of presheaf type if it is classified by a presheaf topos. The class of theories of presheaf type contains all the cartesian theories, as well as many other significant theories (we refer the reader to section \ref{presheaf} for a more extensive discussion on the subject of theories of presheaf type). The classifying topos of a theory of presheaf type $\mathbb T$ is given by the functor category $[\textrm{f.p.} {\mathbb T}\textrm{-mod}(\Set), \Set]$, where $\textrm{f.p.} {\mathbb T}\textrm{-mod}(\Set)$ is the category of (representatives of isomorphism classes of) finitely presentable models of $\mathbb T$. 

Now, if a geometric theory is a quotient (i.e. a theory obtained by adding geometric sequents over the same signature) of a theory of presheaf type $\mathbb T$ then its classifying topos is a subtopos $\Sh(\textrm{f.p.} {\mathbb T}\textrm{-mod}(\Set)^{\textrm{op}}, J)$ of the classifying topos $[\textrm{f.p.} {\mathbb T}\textrm{-mod}(\Set), \Set]$ of $\mathbb T$, and the Grothendieck topology $J$ can be calculated directly by rearranging the axioms of the theory in a particular form involving the formulae which present a (finitely presentable) model of $\mathbb T$ (cfr. Chapter 5 of \cite{OC} for more details). 

Compared with the above-mentioned method via syntactic sites, this latter method of construction of classifying toposes is semantical in spirit, even though, for any theory of presheaf type $\mathbb T$, the category $\textrm{f.p.} {\mathbb T}\textrm{-mod}(\Set)$ can be identified with the opposite of a full subcategory of the geometric syntactic category of $\mathbb T$ (cfr. Theorem 10.3.3 \cite{OC} or section \ref{presheaf} below).

Of course, there are many other methods to prove that a certain\\ Grothendieck topos classifies a given geometric theory; sometimes, one can directly prove that the models of the geometric theory can be identified (naturally in any Grothendieck topos) with flat $J$-continuous functors on a small category $\cal C$, in which case one conclude the classifying topos of the theory is the topos $\Sh({\cal C}, J)$ of sheaves on the site $({\cal C}, J)$. Also, by working with toposes, one is often able to prove the equivalence of a topos presented in terms of one site with a topos presented by using a different site; and, since the notion of classifying topos is clearly invariant under categorical equivalence, this can lead to many different `presentations' of the classifying topos of a given theory.                                             
  
Another source of different sites of definition for the classifying topos of a geometric theory comes from the consideration of quotients of geometric theories. By the duality theorem (cfr. section \ref{duality} below), quotients of a geometric theory $\mathbb T$ can be identified with the subtoposes of the classifying topos of $\mathbb T$. Now, the notion of subtopos is a topos-theoretic invariant (i.e. it depends only on the topos and not on the particular sites of definitions of it) which behaves quite naturally with respect to sites (indeed, for any site $({\cal C}, J)$, the subtoposes of the topos $\Sh({\cal C}, J)$ correspond bijectively with the Grothendieck topologies on $\cal C$ which contain $J$); hence, any particular representation of the classifying topos of a geometric theory yields a related representation of the classifying topos of any quotient of it. 
  
We can think of each site of definition of the classifying topos of a geometric theory as representing a particular aspect of the theory, and of the classifying topos as embodying those essential features of the theory which are invariant with respect to particular (syntactic) presentations of the theory which induce Morita-equivalences at the semantical level. We shall come back to this point in sections \ref{onetopos} and \ref{bridge} below.   
   
We emphasize that there can be really very many different sites of definition for a given Grothendieck topos (the logical interpretation of this fact is central for our purposes and will be discussed in section \ref{morita}). Often, just looking at the same theory in two different ways leads to two different representations of its classifying topos. For example, the coherent theory of fields can be clearly presented both as a quotient of the theory of commutative rings with unit and as a quotient of the theory of von Neumann regular rings, and each of these presentations produces a different representation of the classifying topos (as a category of sheaves on the opposite of the category of finitely presented rings and as a category of sheaves on the opposite of the category of finitely presented von Neumann regular rings). This technical `flexibility' of the theory of toposes in accommodating and extracting the mathematical substance of the apparently immaterial experience of `looking at the same thing in two - or more - different ways' is one of the most striking aspects of the theory. As we shall argue below, an incredible amount of information relevant for classical mathematics is `hidden' inside toposes and can be extracted by using their different sites of definition.

\section{Morita-equivalences}\label{morita}

We have seen in section \ref{logGrothendieck} that two geometric theories have equivalent classifying toposes if and only if they are Morita-equivalent.
Notice that the relation `to be Morita-equivalent to each other' defines an equivalence relation on the collection of all geometric theories, and Grothendieck toposes can be taken as \emph{canonical representatives} of the resulting equivalence classes. Theories which are Morita-equivalent to each other are, broadly speaking, theories which, albeit possibly having a different linguistic (i.e. syntactic) presentation, share a common `semantical core', this core being precisely embodied by their common classifying topos.

Morita-equivalence is a general notion of equivalence of mathematical theories which is ubiquitous in Mathematics (even though there has not been much interest in the past in identifying Morita-equivalences possibly due to the lack of a general theory ascribing central importance to this notion and demonstrating its technical usefulness - one of the purposes of the present paper is in fact to advocate the extreme importance of investigations in this area, cfr. sections \ref{bridge} and \ref{work}). As a simple example of theories which are Morita-equivalent, one can take the theory of Boolean algebras and the theory of Boolean rings.

From the point of view of the logician, it is quite natural to wish to regard as equivalent two mathematical theories whose (categories of) set-based models can be identified with each other; the difference between this notion of equivalence and the notion of Morita-equivalence is simply that in the latter case we require the identification of the models of the two theories to `carry over' to any Grothendieck topos $\cal E$, naturally in $\cal E$. This might seem at first sight a very severe constraint; but in fact, we can expect most of the equivalences of the first kind which arise in mathematical practice to be extensible to Morita-equivalences. The reason for this is that, in establishing an equivalence of the first kind, one generally uses standard set-theoretic constructions which do not involve the law of excluded middle, and, since a Grothendieck topos behaves logically as a `generalized universe of sets' in which one can perform most of the classical set-theoretic constructions with the only significant exception of arguments requiring the law of excluded middle, we can naturally expect to be able to `lift' such equivalence to an arbitrary Grothendieck topos in such a way to globally obtain a Morita-equivalence.    
   
What may perhaps be a bit surprising is that theories of continuous flat functors indeed play a central role into the intrinsically logical subject of Morita-equivalence (cfr. below), also considered that their axiomatizations are fairly bizarre according to the standards of classical model theory (signatures with a possibly infinite number of sorts, infinitary disjunctions etc.). Indeed, $({\cal C}, J)$ is a site of definition for the classifying topos of a geometric theory $\mathbb T$ if and only if $\mathbb T$ is Morita-equivalent to the theory of $J$-continuous flat functors on $\cal C$. In particular, any geometric theory is Morita-equivalent, in a canonical way, to a theory of $J$-continuous flat functors on $\cal C$, namely the theory of $J_{\mathbb T}$-continuous flat functors on ${\cal C}_{\mathbb T}$. 

The notion of Morita-equivalence is also directly connected with that of biinterpretability in classical model theory: in fact, two geometric theories are Morita-equivalent if and only if they are biinterpretable in each other in a generalized sense (a geometric morphism between the classifying toposes of two geometric theories can be regarded as a generalized interpretation of one theory into the other, cfr. section 2.1.5 \cite{OC}). If the theories are coherent this precisely amounts to saying that their syntactic pretoposes are equivalent, which is again a generalized notion of biinterpretability. 

The subject of Morita-equivalence is strictly related to the existence of different sites of definition for a given topos. Indeed, two sites $({\cal C}, J)$ and $({\cal C}', J')$ give rise to the same topos (i.e. $\Sh({\cal C}, J)$ is equivalent to $\Sh({\cal C}', J')$) if and only if the theories of $J$-continuous flat functors on $\cal C$ and of $J'$-continuous flat functors on ${\cal C}'$ are Morita-equivalent; on the other hand, two different mathematical theories which are Morita-equivalent to each other yield different sites of definition of their classifying topos, namely their syntactic sites. In light of this link with the notion of site and of the discussions in section \ref{onetopos}, we can say that the notion of Morita-equivalence indeed captures much of the intuitive idea of `looking at the same thing in different ways'. (Notice that a theory \emph{alone} generates an infinite number of Morita-equivalences: as we remarked in section \ref{onetopos}, just looking at a theory as a quotient of one theory or another leads to a new site of definition for its classifying topos i.e. to a Morita-equivalence.) In fact, a given mathematical property can manifest itself in several different forms in the context of mathematical theories which have a common `semantical core' but a different linguistic presentation; the remarkable fact is that if the property is formulated as a topos-theoretic invariant on some topos then the expression of it in terms of the different theories classified by the topos is determined to a great extent by the technical relationship between the topos and the different sites of definition for it (cfr. section \ref{bridge}).

In connection with this, it seems illuminating to remark that humans have a natural tendency to \emph{visualize} the `semantics' and to use the `syntax' to reason \emph{linguistically} about it; so `looking at the same thing in different ways' may well be taken to mean `describing a given structure by using different languages'. In fact, it often happens in Mathematics (as well as in real life) that the languages and methods used to study some particular object may differ so much from person to person that it could become very difficult to \emph{identify} that in fact the actual object of study is the same. 

To sum up, the fact that different theories are Morita-equivalent to each other translates topos-theoretically in the existence of different sites of definition for one classifying topos. At this point, one might naturally wonder if this link with the notion of site is of any usefulness for the investigation of Morita-equivalences. Actually, the main purpose of this paper is to give a positive answer to this question, on the strength of our new methodologies of topos-theoretic nature for investigating Morita-equivalences. These methodologies, which we will present in the course of the paper, are based on a view of toposes as `bridges' that can be used for transferring information between theories that are Morita-equivalent to each other, and the notion of site plays a central role in this view (cfr. section \ref{bridge} below). In fact, these methods provide a set of techniques for `unifying Mathematics' in the sense of identifying new connections between distinct mathematical theories and \emph{translating} ideas and results between them.

\section{Toposes as `bridges'}\label{bridge}

We have already remarked that there can be many different sites of definition for a given Grothendieck topos, and that this corresponds at the logical level to existence of Morita-equivalences between theories classified by that topos. So, while the assignment of the topos $\Sh({\cal C}, J)$ to a site $({\cal C}, J)$ is a perfectly canonical process, finding a theory classified by a given Grothendieck topos is not at all canonical in general, since it corresponds to finding a small site of definition for the topos: many different sites are sent, via the sheaf construction, to the same topos (up to categorical equivalence). The operation $({\cal C}, J)\rightarrow \Sh({\cal C}, J)$ of taking sheaves on a given site thus appears as a sort of `coding' which extracts exactly those essential features of the theories classified by that topos which are invariant under Morita-equivalence. In a sense, classifying toposes embody the `common features' of geometric theories which are Morita-equivalent to each other.  

In fact, in view of the above-mentioned considerations, it would not even make sense to look for a `privileged' site of definition of a given Grothendieck topos, since this would correspond to a way of canonically selecting a theory out of a Morita-equivalence class, and this is \emph{a priori} an irrational demand, since there is no reason for why one should in general prefer one theory over another (in general, it clearly does not make sense to affirm the `superiority' of one branch of Mathematics over another - all that one can rationally say is that a certain kind of language could be more appropriate than another \emph{in given context}, but these are subjective and contingent considerations).     

Specifically, since the classifying topos of a geometric theory $\mathbb T$ can be taken as a canonical representative for the equivalence class of theories which are Morita-equivalent to $\mathbb T$ (cfr. section \ref{onetopos}), properties of $\mathbb T$ which are invariant with respect to Morita-equivalence are, at least \emph{conceptually}, properties of the classifying topos of $\mathbb T$; conversely, any property of the classifying topos of $\mathbb T$ gives rise to a property of the theories classified by it which is invariant under Morita-equivalence. \emph{Technically}, considered the richness and flexibility of topos-theoretic methods, we can expect these properties of geometric theories which are invariant under Morita-equivalence to be in most cases expressible as invariant properties of their classifying toposes written in topos-theoretic language.
 
Now, the fundamental idea is the following: if we are able to express a property of a given geometric theory as a property of its classifying topos then we can attempt to express this property in terms of any of the other theories having the same classifying topos, so to obtain a relation between the original property and a new property of a different theory which is Morita-equivalent to it. The classifying topos thus acts as a sort of `bridge' connecting different mathematical theories that are Morita-equivalent to each other, which can be used to transfer information and results from one theory to another. The purpose of the present paper is to show that this idea of toposes as unifying spaces is technically very feasible; the great amount and variety of results in the Ph.D. thesis \cite{OC} give clear evidence for the fruitfulness of this point of view and, as we shall argue in the course of the paper, a huge number of new insights into any field of Mathematics can be obtained as a result of the application of these techniques.

Indeed, the fact that different mathematical theories can have equivalent classifying toposes translates into the existence of different sites of definition for one topos. Topos-theoretic invariants (i.e. properties of toposes which are invariant under categorical equivalence) can then be used to transfer properties from one theory to another. As we shall see in section \ref{invar}, this is made possible by the fact that the abstract relationship between a site $({\cal C},J)$ and the topos $\Sh({\cal C}, J)$ which it `generates' is often very natural (in the sense that properties of sites technically relate to topos-theoretic invariants in a natural way), enabling us to easily transfer invariants across different sites. Notice that, since the construction of the category of sheaves $\Sh({\cal C}, J)$ from a site $({\cal C},J)$ is entirely canonical, a property of the topos $\Sh({\cal C}, J)$ \emph{is}, at least in principle, a property of the site $({\cal C},J)$. What happens in practice is that these properties of sites often have a genuine `categorical' description, or at least are implied or imply a property admitting such a description. For example, if $\cal C$ is a (small) category satisfying the right Ore condition and $J$ is the atomic topology on it then the topos $\Sh({\cal C}, J)$ is atomic (notice that atomicity is a topos-theoretic invariant). For many invariants (for example, the property of a topos to be Boolean, De Morgan, or two-valued), one has bijective characterizations of the kind ``$\Sh({\cal C}, J)$ satisfies the invariant if and only if the site $({\cal C}, J)$ satisfies a certain `tractable' categorical property'' for an arbitrary site $({\cal C}, J)$, which allow a direct transfer of information between distinct sites of definition of the same topos. For other invariants, one may have implications going in general just in one direction, while for going in the other direction one has to assume that the site is of a particular kind (for example, subcanonical); as we shall see in section \ref{logic} below, the geometric syntactic sites of geometric theories behave particularly well in relation to such characterizations. Of course, it may also happen be that one cannot establish for a given invariant general `tractable' site characterizations of the above-mentioned kind, but that in the particular cases of interest one can employ \emph{ad hoc} arguments to identify properties of the site which imply or follow from the fact that the topos satisfies the given invariant. We shall see these methodologies in action in a variety of different contexts in the course of the paper.          

As we shall see in sections \ref{invar} and \ref{duality} below, the level of generality represented by topos-theoretic invariants is ideal to capture several important features of mathematical theories. Indeed, topos-theoretic invariants considered on the classifying topos $\Set[{\mathbb T}]$ of a geometric theory $\mathbb T$ translate into interesting logical (i.e. syntactic or semantic) properties of $\mathbb T$, and topos-theoretic constructions on classifying toposes correspond to natural operations on the theories classified by them.   

To sum up, it is by means of expressing a topos-theoretic invariant in terms of the different sites of definition of a given topos that the transfer of information between theories classified by that topos takes place; it is precisely in this sense that toposes act as `bridges' connecting different mathematical theories. The notion of Morita-equivalence thus acquires a prominent role in our context; it is in a sense the primitive ingredient on which the machinery just described can be put at work. It is therefore natural to wonder how one can find Morita-equivalences starting from one's mathematical work. There are indeed many ways through which one can arrive at Morita-equivalences. We have already seen that Topos Theory itself is a primary source of Morita-equivalences, since any alternative way of representing a topos as a category of sheaves on a site leads to a Morita-equivalence (cfr. sections \ref{onetopos} and \ref{morita} above). The logical approach to Morita-equivalences and the problem of extending a classical equivalence to a Morita-equivalence was already discussed in section \ref{morita}. Sheaf representations for various kind of structures should also admit a natural topos-theoretic interpretation as Morita-equivalences. Moreover, it is reasonable to expect that most of the classical \emph{dualities} arising in Mathematics which involve some geometric theory should in some way give rise to Morita-equivalences extracting the essential features of them.

A natural way by means of which the `working mathematician' can enter the world of toposes so to benefit from the existing knowledge on Morita-equivalences and the application of the above-mentioned methods, is the following: whenever he or she finds a mathematical property of his or her interest, he or she should try to find a site $({\cal C}, J)$ and an invariant of the topos $({\cal C}, J)$ which relates (in the sense of implying or being implied, at least under additional assumptions) to the original property. For example, given a small category ${\cal C}$ satisfying the amalgamation property (cfr. section \ref{fraisse} below), the property of ${\cal C}$ to satisfy the joint embedding property (cfr. section \ref{fraisse} below) is easily seen to be equivalent to the property of the topos $\Sh({\cal C}^{\textrm{op}}, J_{at})$, where $J_{at}$ is the atomic topology on ${\cal C}^{\textrm{op}}$, to be two-valued. Since the classifying topos of a geometric theory is two-valued if and only if the theory is complete (in the sense that any geometric sentence over the signature of the theory is provably equivalent to $\top$ or $\bot$, but not both), we conclude that the joint embedding property on our category $\cal C$, which is a property of `\emph{geometrical}' flavour, translates in the \emph{logical} property of completeness of the theory of $J_{at}$-continuous flat functors on ${\cal C}^{\textrm{op}}$. This example shows that common mathematical properties (in this case, the joint embedding property on a category $\cal C$) may well arise as specializations of abstract logical properties of theories which are stable under Morita-equivalence (in this case, the property of completeness of geometric theories) to a particular geometric theory (in this case, the theory of $J_{at}$-continuous flat functors on ${\cal C}^{\textrm{op}}$). We will come back to this example in our discussion of the topos-theoretic interpretation of Fra\"iss\'e's construction in section \ref{fraisse}. Notice the central role played by the theories of continuous flat functors in these issues; they indeed have a strong link with classical mathematics, and part of the reason for this resides in the fact that their axiomatizations directly involve the objects and arrows of a given category, respectively as sorts and function symbols in their signatures. 

In general, inventing arguments as the one just described is not as difficult as it might seem, and it will become easier and easier with the time, as topos-theorists will discover new invariants for the benefit of mathematicians. Anyway, it is the author's opinion that this should be a participative effort; as topos-theorists have drawn a lot of inspiration from general topology in designing the topological properties of toposes, so in the future they could be motivated in their work by the needs of people working in any mathematical field. As we shall see in section \ref{invar} below, the existence of topos-theoretic invariants with particular properties can have important ramifications in specific mathematical contexts; toposes indeed occupy a central role in Mathematics, and their use can be strategic in several situations.  

The application of the methodologies just described enables one to extract an incredible amount of new information on Morita-equivalences, and to establish connections between distinct theories which could hardly be visible otherwise. In fact, the kind of insight that these methods can generate is intrinsically different from, if not actually subsuming (it may a bit early for saying this), that provided by the traditional methods of transferring information between equivalent theories by using the specific description of the equivalence. Indeed, in our approach, instead of using the explicit description of the Morita-equivalence, one exploits the relationship between the topos and its sites of definition. This is possible because we have an actual mathematical object, namely the classifying topos, which extracts the `common core' of the theories classified by it, and the equivalences between the theories are essentially coded in the relationship between the topos and its different sites of definition. In fact, for most purposes, it is only the \emph{existence} of a Morita-equivalence that really matters, and we can well ignore the actual description of it; of course, if one wants to establish more `specific' results, an explicit description of the Morita-equivalence becomes necessary (in which case, one can use invariant properties of \emph{objects} of toposes rather than invariants properties of the `whole topos', cfr. section \ref{examp} below) but for treating most of the `global' properties of theories this is not at all necessary (since, by definition, a topos-theoretic invariant is stable under \emph{any} kind of categorical equivalence). In other words, we can generate a lot of interesting insights still remaining at the `$1$-dimensional level'. We emphasize that there is an strong element of \emph{automatism} in the techniques just explained; by means of these methods, one can generate new mathematical results without really making any conscious effort: indeed, in most cases one can just readily apply the well-known general characterizations connecting properties of sites and topos-theoretic invariants (such as for example the above-mentioned ones for the property of Boolean and two-valuedness) to the particular case of interest. On the other hand, the range of applicability of these methods is boundless within Mathematics, by the very generality of the notion of topos.     

The investigation of theories of presheaf type carried out in the author's Ph.D. thesis and briefly discussed in section \ref{presheaf} below, represents a clear illustration of these principles. In fact, if a geometric theory $\mathbb T$ is of presheaf type, we automatically have two different representations of its classifying topos: $[\textrm{f.p.} {\mathbb T}\textrm{-mod}(\Set), \Set]$ and $\Sh({\cal C}_{\mathbb T}, J_{\mathbb T})$. And, starting from this dual representation of the classifying topos, one can obtain a great number of insights by applying these methodologies, including a general version of Fra\"iss\'e's theorem in classical model theory (cfr. section \ref{fraisse}).

\section{Topos-theoretic invariants}\label{invar}

We have already explained the role of topos-theoretic invariants in our view of toposes as unifying spaces which can act as `bridges' connecting distinct mathematical theories which are Morita-equivalent to each other. But, what do we exactly mean by the expression `topos-theoretic invariant'? By this expression, we mean any kind of \emph{property} $P$ of (families of) Grothendieck toposes or \emph{construction} $C$ involving (families of) Grothendieck toposes which is invariant under categorical equivalence of Grothendieck toposes i.e. such that if a given family of Grothendieck toposes $\{{\cal E}_{i} \textrm{ | } i \in I\}$ satisfies the property $P$ then for any family $\{{\cal E}_{i}' \textrm{ | } i \in I\}$ of toposes such that for each $i\in I$ the topos ${\cal E}_{i}$ is equivalent to the topos ${\cal E}_{i}'$, the latter family also satisfies the property $P$ (the requirement for a construction $C$ is that if $\{{\cal E}_{i} \textrm{ | } i \in I\}$ and $\{{\cal E}_{i}' \textrm{ | } i \in I\}$ are two families of toposes such that for each $i\in I$ the topos ${\cal E}_{i}$ is equivalent to the topos ${\cal E}_{i}'$ then the result of applying the construction $C$ to the first family should be equivalent - in a mathematical sense - to the result of applying $C$ to the second family). We do \emph{not} require that there should be a linguistic expression of the property $P$ (resp. a description of the construction $C$) in the language of Topos Theory, all that matters is that $P$ should be invariant under categorical equivalence of toposes. Of course, any property which admits a linguistic description in the (informal) language of Topos Theory is automatically a topos-theoretic invariant (as a consequence of a very general `meta-theorem'), but one should not be limited by this thought in looking for invariants. In fact, we can expect there to be many more topos-theoretic invariants around than the ones that we can presently describe in the categorical language that we commonly use for talking and reasoning about toposes.  

Examples of well-known topos-theoretic invariants include: the property of a topos to be \emph{Boolean}, \emph{De Morgan}, \emph{atomic}, \emph{two-valued}, to be \emph{connected}, to be \emph{locally connected}, to be \emph{compact}, to be \emph{local}, to be a \emph{subtopos} of a given topos (in the sense of geometric inclusion), to the the \emph{classifying topos} of a given geometric theory, to have enough points, etc. Cohomology and homotopy groups of toposes are also important topos-theoretic invariants. 

As we have anticipated in the last section, most of these invariants behave quite naturally with respect to sites, that is we have general characterizations connecting `natural' properties of sites with these invariants on the corresponding toposes of sheaves. In fact, for most of the `topologically motivated' invariants of toposes, we have natural characterizations going in one direction i.e. asserting implications of the kind ``if a site $({\cal C}, J)$ satisfies a given property then the topos $\Sh({\cal C}, J)$ satisfies the given invariant'' (we refer to part C of \cite{El} for a detailed presentation of these results), and, under specific assumptions on the site, one can often establish results going in the converse direction (cfr. section \ref{logic} below). On the other hand, for other `logically motivated' invariants such as the property of a topos to be Boolean, De Morgan, or two-valued, we have natural \emph{bijective} characterizations (cfr. \cite{OC}, in particular chapters $6$ and $9$). In section \ref{duality}, we will focus on a particular invariant, namely the notion of \emph{subtopos}, and discuss its relevance for classical Mathematics.

As we have already remarked in the last section, it is often useful to consider invariants of objects of toposes, in order to establish `local' properties rather than `global' ones. By an `invariant of objects of toposes' we mean a property $Q$ of objects of toposes such that whenever $\tau:{\cal E}\to {\cal F}$ is an equivalence of toposes and $a$ is an object of ${\cal E}$, the object $a$ satisfies $Q$ if and only if the object $\tau(a)$ satisfies $Q$ (and similarly for families of objects). Examples of invariants of object of toposes include: the property of an object to be compact, connected, indecomposable, irreducible, coherent, to be an atom, etc. (cfr. chapter 10 of \cite{OC} for a formal definition of these invariants). 

We have seen that any topos-theoretic invariant can generate, via the methodologies described above, a considerable amount of mathematical results. On the other hand, also `negative results' asserting the \emph{non-existence} of topos-theoretic invariants with particular properties can be directly connected to relevant results in Mathematics. For example, one might wonder whether there exists a topos-theoretic invariant $P$ with the following property: whenever $\cal E$ is the topos of sheaves $\Sh_{G}(X)$ on a topological groupoid $s,t:G\to X$, $\cal E$ satisfies $P$ if and only if the groupoid $s,t:G\to X$ is `algebraically connected' i.e. the coequalizer of $s$ and $t$ in the category of topological spaces is given by the one-point space. An invariant satisfying the analogous of this property for localic groupoids exists, and in fact it is the precisely the property of a topos to be hyperconnected (cfr. Lemma C5.3.7 \cite{El}). On the other hand, no such invariant can exist for topological groupoids, and the non-existence of such invariant is directly related to the classical model-theoretic fact that there can be (coherent) theories which are \emph{not} countably categorical but which are $k$-categorical for some uncountable cardinal $k$ (or viceversa). Indeed, from the representation theorem by I. Moerdijk and C. Butz of Grothendieck toposes with enough points as toposes of sheaves on a topological groupoid (cfr. \cite{BM}), it follows that the classifying topos of a coherent theory $\cal T$ can be represented as the topos of sheaves $\Sh_{G}(X)$ on the topological groupoid $s,t:G\to X$ of isomorphisms of models of $\mathbb T$ in $X$, where $X$ is \emph{any} set of (enumerated) models of $\mathbb T$ in $\Set$ which are jointly conservative for $\mathbb T$. Now, since by the L\"owenheim-Skolem theorem the collection of all models of a countable theory of a fixed cardinality (either countably infinite or uncountable) is jointly conservative for the theory, the existence of our invariant $P$ would imply that for any countable coherent theory $\mathbb T$, $\mathbb T$ has exactly one isomorphism class of models of cardinality $k$ if and only if it has exactly one isomorphism class of models of cardinality $k'$, for \emph{any} cardinals $k,k' \geq \omega$, and this is well-known not to be true in general (a counterexample is given by the coherent theory of algebraically closed fields of a given characteristic).

\subsection{The logical meaning of invariants}\label{logic}

In view of our previous claims concerning the importance of topos-theoretic invariants for the investigation of geometric theories, it is natural to wonder whether important invariants on classifying toposes relate to interesting logical (i.e. syntactic or semantic) properties of the theories classified by them. The purpose of the previous section is to show that they do. In fact, this line of investigation has been systematically carried out in \cite{OC}, mostly in chapters 6 and 10, where many well-known invariants were studied from the point of view of geometric theories. The main tool used in these investigations is the geometric syntactic site of a geometric theory. Recall that the classifying topos of a geometric theory $\mathbb T$ can always be represented as the category $\Sh({\cal C}_{\mathbb T}, J_{\mathbb T})$ of sheaves on the geometric syntactic site $({\cal C}_{\mathbb T}, J_{\mathbb T})$ of $\mathbb T$. Now, this site behaves quite well in relation to the problem of transferring properties from the topos to the site. Indeed, the Grothendieck topology $J_{\mathbb T}$ is subcanonical and the Yoneda embedding $y:{\cal C}_{\mathbb T} \to \Sh({\cal C}_{\mathbb T}, J_{\mathbb T})$ satisfies the additional useful property that any subobject in $\Sh({\cal C}_{\mathbb T}, J_{\mathbb T})$ of an object of the form $y(c)$ is of the form $y(d)\mono y(c)$ for some subobject $d\mono c$ in $\Sh({\cal C}_{\mathbb T}, J_{\mathbb T})$. Thus, since properties of the syntactic site of a theory often rephrase naturally as syntactic properties of the theory, the process of expressing topos-theoretic invariants on the classifying topos of a theory $\mathbb T$ as syntactic properties of $\mathbb T$ works in general quite smoothly. In order to illustrate these points, we report below some of the results from \cite{OC} concerning the logical interpretation of invariants. We emphasize that these results hold \emph{uniformly} for any geometric theory i.e. no particular properties of the theory are assumed in order to achieve the characterizations. In order to state the results, we first have to introduce some terminology.

\begin{definition}
Let $\mathbb T$ be a geometric theory over a signature $\Sigma$ and $\phi(\vec{x})$ a geometric formula-in-context over $\Sigma$. Then

\begin{enumerate}[(i)]

\item We say that $\phi(\vec{x})$ is \emph{$\mathbb T$-complete} if the sequent ($\phi \vdash_{\vec{x}} \bot$) is not provable in $\mathbb T$, and for every geometric formula $\chi(\vec{x})$ in the same context either ($\phi \vdash_{\vec{x}} \chi$) or ($\chi \wedge \phi \vdash_{\vec{x}} \bot$) is provable in $\mathbb T$.

\item We say that $\phi(\vec{x})$ is \emph{$\mathbb T$-indecomposable} if for any family $\{\psi_{i}(\vec{x}) \textrm{ | } i\in I\}$ of geometric formulae in the same context such that for each $i$, $\psi_{i}$ $\mathbb T$-provably implies $\phi$ and for any distinct $i, j\in I$, $\psi_{i}\wedge \psi_{j} \vdash_{\vec{x}} \bot$ is provable in $\mathbb T$, we have that $\phi \vdash_{\vec{x}} \mathbin{\mathop{\textrm{\huge $\vee$}}\limits_{i\in I}} \psi_{i}$ provable in $\mathbb T$ implies $\phi \vdash_{\vec{x}} \psi_{i}$ provable in $\mathbb T$ for some $i\in I$.

\item We say that $\phi(\vec{x})$ is \emph{$\mathbb T$-irreducible} if for any family $\{\theta_{i} \textrm{ | } i\in I\}$ of $\mathbb T$-provably functional geometric formulae $\{\vec{x_{i}}, \vec{x}.\theta_{i}\}$ from $\{\vec{x_{i}}. \phi_{i}\}$ to $\{\vec{x}. \phi\}$ such that $\phi \vdash_{\vec{x}} \mathbin{\mathop{\textrm{\huge $\vee$}}\limits_{i\in I}}(\exists \vec{x_{i}})\theta_{i}$ is provable in $\mathbb T$, there exist $i\in I$ and a $\mathbb T$-provably functional geometric formula $\{\vec{x}, \vec{x_{i}}. \theta'\}$ from $\{\vec{x}. \phi\}$ to $\{\vec{x_{i}}. \phi_{i}\}$ such that $\phi \vdash_{\vec{x}} (\exists \vec{x_{i}})(\theta' \wedge \theta_{i})$ is provable in $\mathbb T$.

\item We say that $\phi(\vec{x})$ is \emph{$\mathbb T$-compact} if for any family $\{\psi_{i}(\vec{x}) \textrm{ | } i\in I\}$ of geometric formulae in the same context, $\phi \vdash_{\vec{x}} \mathbin{\mathop{\textrm{\huge $\vee$}}\limits_{i\in I}} \psi_{i}$ provable in $\mathbb T$ implies $\phi \vdash_{\vec{x}} \mathbin{\mathop{\textrm{\huge $\vee$}}\limits_{i\in I'}} \psi_{i}$ provable in $\mathbb T$ for some finite subset $I'$ of $I$.

\end{enumerate}

\end{definition}

These notions arise precisely from the rephrasing of invariant properties of objects of toposes, namely the property on an object to be an atom, to be indecomposable, to be irreducible, to be compact, in terms of objects of the syntactic site of the theory; as such, they represent an implementation of the abstract methodologies presented in section \ref{bridge} in the case of the classifying topos of a geometric theory represented in terms of the syntactic site of the theory. Specifically, if $y:{\cal C}_{\mathbb T} \to \Sh({\cal C}_{\mathbb T}, J_{\mathbb T})$ is the Yoneda embedding, a geometric formula-in-context $\{\vec{x}. \phi\}$ is $\mathbb T$-complete (resp. $\mathbb T$-indecomposable, $\mathbb T$-irreducible, $\mathbb T$-compact) if and only if the object $y(\{\vec{x}. \phi\})$ in the topos $\Sh({\cal C}_{\mathbb T}, J_{\mathbb T})$ is an atom (resp. indecomposable, irreducible, compact). The last step to get the link with the invariants on the classifying topos of $\mathbb T$ is to relate these `global' properties of toposes to the above-mentioned `local' properties of objects; this turns out to be possible in these cases, and the outcome is the following.    

\begin{theorem}
Let $\mathbb T$ be a geometric theory over a signature $\Sigma$ and $\Set[{\mathbb T}]$ be its classifying topos. Then

\begin{enumerate}

\item $\Set[{\mathbb T}]$ is \emph{locally connected} (resp. \emph{atomic}) if and only if for any geometric formula $\phi(\vec{x})$ over $\Sigma$ there exists a (unique) family $\{\psi_{i}(\vec{x}) \textrm{ | } i\in I\}$ of $\mathbb T$-indecomposable (resp. $\mathbb T$-complete) geometric formulae in the same context such that

\begin{enumerate}[(i)]
\item for each $i$, $\psi_{i}$ $\mathbb T$-provably implies $\phi$,

\item for any distinct $i, j\in I$, $\psi_{i}\wedge \psi_{j} \vdash_{\vec{x}} \bot$ is provable in $\mathbb T$, and

\item $\phi \vdash_{\vec{x}} \mathbin{\mathop{\textrm{\huge $\vee$}}\limits_{i\in I}} \psi_{i}$ is provable in $\mathbb T$. 
\end{enumerate}

\item $\Set[{\mathbb T}]$ is \emph{equivalent to a presheaf topos} if and only if there exists a collection $\cal F$ of $\mathbb T$-irreducible geometric formulae-in-context over $\Sigma$ satisfying the following property: for any geometric formula $\{\vec{y}. \psi\}$ over $\Sigma$, there exist objects $\{\vec{x_{i}}. \phi_{i}\}$ in $\cal F$ as $i$ varies in $I$ and $\mathbb T$-provably functional geometric formulae $\{\vec{x_{i}}, \vec{y}.\theta_{i}\}$ from $\{\vec{x_{i}}. \phi_{i}\}$ to $\{\vec{y}. \psi\}$ such that $\psi \vdash_{\vec{y}} \mathbin{\mathop{\textrm{\huge $\vee$}}\limits_{i\in I}}(\exists \vec{x_{i}})\theta_{i}$ is provable in $\mathbb T$;

\item $\Set[{\mathbb T}]$ is \emph{compact} if and only if the formula $\top$ in the empty context is $\mathbb T$-compact. 

\item $\Set[{\mathbb T}]$ is \emph{two-valued} if and only if the formula $\top$ in the empty context is $\mathbb T$-complete. 

\end{enumerate}

\end{theorem}
  
Other results concerning the logical interpretation of topos-theoretic invariants which were obtained in \cite{OC} are the following:

\begin{theorem}
Let $\mathbb T$ be a geometric theory over a signature $\Sigma$ and $\Set[{\mathbb T}]$ be its classifying topos. Then

\begin{enumerate}

\item $\Set[{\mathbb T}]$ is \emph{Boolean} if and only if for every geometric formula $\phi(\vec{x})$ over $\Sigma$ there is a geometric formula $\chi(\vec{x})$ over $\Sigma$ in the same context such that $\phi(\vec{x}) \wedge \chi(\vec{x}) \:\vdash_{\vec{x}}\: \bot$ and $\top \:\vdash_{\vec{x}}\: \phi(\vec{x}) \vee \chi(\vec{x})$ are provable in $\mathbb T$.  

\item $\Set[{\mathbb T}]$ is \emph{De Morgan} if and only if for every geometric formula $\phi(\vec{x})$ over $\Sigma$, there exist two consistent geometric formulae $\psi_{1}(\vec{x})$ and $\psi_{2}(\vec{x})$ over $\Sigma$ in the same context such that:\\
$\top \vdash_{\vec{x}} \psi_{1}(\vec{x})\vee \psi_{2}(\vec{x})$ is provable in $\mathbb T$,\\
$\psi_{1}(\vec{x})\wedge \phi(\vec{x}) \vdash_{\vec{x}} \bot$ is provable in $\mathbb T$ and\\
for every geometric formula $\chi(\vec{x})$ over $\Sigma$ in the same context such that $\chi(\vec{x}) \vdash_{\vec{x}} \psi_{2}(\vec{x})$ is provable in $\mathbb T$, $\chi(\vec{x}) \wedge \phi(\vec{x}) \vdash_{\vec{x}} \bot$ is provable in $\mathbb T$ if and only if $\chi(\vec{x}) \vdash_{\vec{x}} \bot$ is provable in $\mathbb T$.  

\item $\Set[{\mathbb T}]$ has \emph{enough points} if and only if whenever a geometric sequent $\sigma$ over $\Sigma$ is valid in every $\Set$-model of $\mathbb T$, $\sigma$ is provable in $\mathbb T$.

\end{enumerate}

\end{theorem}

These results show that topos-theoretic invariants on the classifying topos of a geometric theory indeed correspond to natural and interesting logical properties of the theory, `uniformly' for any geometric theory; we shall see some applications of these notions in the next sections of the paper.

In connection with this, we also remark that, since every Grothendieck topos has a syntactic site of definition (being the classifying topos of a geometric theory), all the invariant properties or constructions on toposes, and in particular all the results in Topos Theory, have, at least in principle, a logical counterpart; and, since the transfer of properties from the topos to the syntactic site and then to the geometric theory works in general quite smoothly, we can expect most invariants of Grothendieck toposes to give rise to properties of geometric theories expressed in the language of geometric logic. As an example, we mention Deligne's theorem on coherent toposes, which is well-known to be equivalent to the classical completeness theorem for coherent logic (in fact, as we shall see in section \ref{examp}, this equivalence falls under our general scheme for transferring knowledge between distinct mathematical theories by using toposes).  

So (geometric) Logic provides a language in which most of Topos Theory, and hence of Mathematics, can be read and - we emphasize - not only \emph{formulated} but actually \emph{done}. This indicates that Logic in general, and geometric logic in particular, can be fundamental tools for solving a huge variety of mathematical problems. With their work, Model Theorists have already shown the fundamental impact that the logical investigations can have for classical Mathematics; we believe that this positive trend will continue in the future and will be greatly enriched by the contact with the topos-theoretic methodologies.

\newpage 
\section{The duality theorem}\label{duality}

In this section, we discuss the importance of the notion of \emph{subtopos} for geometric logic and Mathematics in general; as we shall see below, this is a topos-theoretic invariant that behaves particularly well with respect to sites; indeed, the subtoposes of a topos $\Sh({\cal C}, J)$ correspond bijectively to the Grothendieck topologies on $\cal C$ which contain $J$. This enables us to fruitfully apply our philosophy `toposes as bridges' in the context of subtoposes.  

In \cite{OC}, a \emph{duality theorem} connecting subtoposes of the classifying topos of a geometric theory and $\mathbb T$ and geometric `quotients' of $\mathbb T$ is proved.  

Before we can state the theorem, we need to introduce a couple of definitions. 

\begin{definition}
Let $\mathbb{T}$ be a geometric theory over a signature $\Sigma$. A \emph{quotient} of $\mathbb{T}$ is a geometric theory $\mathbb{T}'$ over $\Sigma$ such that every axiom of $\mathbb{T}$ is provable in $\mathbb{T}'$.
\end{definition}
 
\begin{definition}
Let $\mathbb{T}$ and $\mathbb{T}$ be geometric theories over a signature $\Sigma$. We say that $\mathbb{T}$ and $\mathbb{T}'$ are \emph{syntactically equivalent}, and we write $\mathbb{T} \equiv_{s} \mathbb{T}'$, if for every geometric sequent $\sigma$ over $\Sigma$, $\sigma$ is provable in $\mathbb{T}$ if and only if $\sigma$ is provable in $\mathbb{T}'$.
\end{definition}      

\begin{theorem}\label{LatticesTheories_dualita}
Let $\mathbb{T}$ be a geometric theory over a signature $\Sigma$. Then the assignment sending a quotient of $\mathbb{T}$ to its classifying topos defines a bijection between the $\equiv_{s}$-equivalence classes of quotients of $\mathbb T$ and the subtoposes of the classifying topos $\Set[\mathbb{T}]$ of $\mathbb{T}$.   
\end{theorem}

The bijection given by the theorem is natural in the following sense. If $i_{{\cal F}}:{\cal F}\hookrightarrow \Set[\mathbb{T}]$ is the subtopos of $\Set[\mathbb{T}]$ corresponding to a quotient ${\mathbb T}'$ of $\mathbb T$ via the duality theorem, we have a commutative (up to natural isomorphism) diagram in $\textbf{CAT}$ (where $i$ is the obvious inclusion)

\[  
\xymatrix {
{{\mathbb T}'}\textrm{-mod}({\cal E}) \ar[rr]^{\simeq} \ar[d]^{i} & & {\bf Geom}({\cal E}, {\cal F}) \ar[d]^{i_{{\cal F}}\circ -} \\
{\mathbb T}\textrm{-mod}({\cal E}) \ar[rr]^{\simeq} & &  {\bf Geom}({\cal E},\Set[\mathbb{T}])}
\]

naturally in ${\cal E}\in \mathfrak{BTop}$.

Notice that an immediate consequence of the duality theorem is the fact that two geometric theories over the same signature have equivalent classifying toposes if and only if they are syntactically equivalent.

In \cite{OC} two different proofs of this theorem were given, one relying on the theory of classifying toposes and the other one based on a proof-theoretic interpretation of the notion of Grothendieck topology (which we shall discuss below). Anyway, both arguments exploit in an essential way the syntactic representation of the classifying topos of $\mathbb T$ as the category of sheaves on the geometric syntactic site of $\mathbb T$. So this result falls under the general pattern discussed in section \ref{logic} of rephrasing topos-theoretic invariants on the classifying topos of a theory $\mathbb T$ in terms of logical properties or constructions on $\mathbb T$ by passing through the geometric syntactic site of $\mathbb T$. 

The duality theorem realizes a unification of the theory of elementary toposes with geometric logic, passing through the theory of Grothendieck toposes. Indeed, the theorem allows us to interpret many concepts of elementary topos theory which apply to the lattice of subtoposes of a given topos at the level of geometric theories. These notions include for example the coHeyting algebra structure on the lattice of subtoposes of a given topos, open, closed, quasi-closed subtoposes, the dense-closed factorization of a geometric inclusion, coherent subtoposes, subtoposes with enough points, the surjection-inclusion factorization of a geometric morphism, skeletal inclusions, atoms in the lattice of subtoposes of a given topos, \emph{Booleanization} and \emph{DeMorganization} of a topos. The remarkable fact is that the resulting notions and results in geometric logic are of considerable logical interest. We shall survey some of them below.

First of all, we remark that the duality theorem has an elegant proof-theoretic interpretation. In order to describe it, we first need to note that the notion of Grothendieck topology on $\cal C$ gives naturally rise to an abstract kind of `proof system'. Specifically, given a collection $\cal A$ of sieves on a given category $\cal C$, we can define a proof system ${\cal T}_{\cal C}^{\cal A}$ as follows: the axioms of ${\cal T}_{\cal C}^{\cal A}$ are the sieves in $\cal A$ together with all the maximal sieves, while the inference rules of ${\cal T}_{\cal C}^{\cal A}$ are the proof-theoretic versions of the well-known axioms for Grothendieck topologies, i.e. the rules:

\emph{Stability rule:}
\[
\begin{array}{c}
R\\
\hline
f^{\ast}(R)  
\end{array}   
\] 
where $R$ is any sieve on an object $c$ in ${\cal C}$ and $f$ is any arrow in ${\cal C}$ with codomain $c$.

\emph{Transitivity rule:}
\[
\begin{array}{c}
Z \textrm{     } \{f^{\ast}(R) \textrm{ | } f\in Z \}\\
\hline
R  
\end{array}   
\] 
where $R$ and $Z$ are sieves in ${\cal C}$ on a given object of $\cal C$.

In these terms, the duality theorem yields, for any geometric theory $\mathbb T$, a kind of \emph{proof-theoretic equivalence} between the traditional proof system of geometric logic with the addition of the axioms of $\mathbb T$, and the system ${\cal T}_{\cal C}^{\cal A}$, where ${\cal C}$ is the geometric syntactic category ${\cal C}_{\mathbb T}$ of the theory $\mathbb T$ and ${\cal A}$ is the collection of sieves in the syntactic topology $J_{\mathbb T}$ of $\mathbb T$. Specifically, one can associate to each geometric sequent over the signature of $\mathbb T$ a sieve in the category ${\cal C}$ and conversely, in such a way that these two assignments are natural with respect to the notions of provability of the two proof systems and are inverse to each other up to provability. The `closed' theories of the first proof system can be clearly identified with the $\equiv_{s}$-equivalence classes of quotients of $\mathbb T$, while the `closed' theories of the second proof system are the Grothendieck topologies on ${\cal C}_{\mathbb T}$ which contain $J_{\mathbb T}$ (for more details, we refer the reader to Chapter 2 of \cite{OC}).

Similar equivalences between classical proof systems of geometric logic and systems associated to Grothendieck topologies arise in the context of theories of presheaf type. In fact, if $\mathbb T$ is a theory of presheaf type then the duality theorem gives rise to a bijection between the ($\equiv_{s}$-equivalence classes of) quotients of $\mathbb T$ and the Grothendieck topologies $J$ on the opposite of the category $\textrm{f.p.} {\mathbb T}\textrm{-mod}(\Set)$ of finitely presentable $\mathbb T$-models. It turns out that this bijection can be made into a proof-theoretic equivalence similarly as above. This fact has important consequences; for example, it follows that any quotient of a theory of presheaf type $\mathbb T$ has a presentation over its signature in which all axioms are of the form $\phi \vdash_{\vec{x}} \mathbin{\mathop{\textrm{\huge $\vee$}}\limits_{i\in I}}(\exists \vec{y_{i}})\theta_{i}$, where, for any $i\in I$, $\theta_{i}(\vec{y_{i}}, \vec{x})$ is a $\mathbb T$-provably functional formula from $\{\vec{y_{i}}. \psi\}$ to $\{\vec{x}. \phi\}$ and $\phi(\vec{x})$, $\psi(\vec{y_{i}})$ are formulae which present a $\mathbb T$-model. 

These `proof-theoretic equivalences' are also very useful in practice because they enable one to work with Grothendieck topologies in order to obtain axiomatizations of geometric theories or syntactic results about them, and Grothendieck topologies have in many contexts a considerable `computational advantage' over the classical proof systems. For example, there is a useful formula for the Grothendieck topology generated by a given family of sieves (this has for example been used in \cite{OC} to establish a deduction theorem for geometric logic) and of the Heyting operation between Grothendieck topologies (cfr. Chapter 3 of \cite{OC} for all these results). An application of this discourse to the calculation of the (lattice-theoretic) meet of geometric theories was given in \cite{OC}, where a natural axiomatization of the meet of the theory of local rings and integral domains is achieved by computing the intersection of the corresponding Grothendieck topologies on the opposite of the category of finitely presented commutative rings with unit (notice that this kind of axiomatization problems are in general far from being trivial). In fact, this example represents an application of our philosophy `toposes as bridges' described in section \ref{bridge}; indeed, in this case the invariant is the meet of two subtoposes, which, looked from the point of view of the syntactic site of the theory $\mathbb R$ of commutative rings with unit amounts to the meet of the two quotients of the theory $\mathbb R$, while looked from the point of view of the trivial site of the classifying topos of $\mathbb R$ (the opposite $\cal R$ of the category of finitely presented rings equipped with the trivial Grothendieck topology) corresponds to the intersection of the associated Grothendieck topologies on the category $\cal R$.  

Let us proceed to describe other applications of the duality theorem. It follows from the theorem that the collection of ($\equiv_{s}$-equivalence classes of) geometric theories over a given signature has the structure of a Heyting algebra with respect to the natural ordering between theories given by `${\mathbb T}'\leq {\mathbb T}''$ if and only if all the axioms of ${\mathbb T}'$ are provable in ${\mathbb T}''$'. Also, it is possible to derive from the theorem explicit descriptions of the lattice and Heyting operations on geometric theories (cfr. \cite{OC}), by using the explicit description of the Heyting operations on Grothendieck topologies. This is another illustration of the fact that working with Grothendieck topologies can be easier than arguing in the classical proof system of geometric logic; indeed, it would be rather intricate to achieve these characterizations by means of the traditional logical methods. 

In light of the fact that the notion of subtopos is a topos-theoretic invariant, the duality theorem allows us to easily transfer information between quotients of geometric theories classified by the same topos. For example, consider the following problem. Let ${\mathbb T}_{1}$ and ${\mathbb T}_{2}$ be two geometric theories equipped with a Morita-equivalence between them; is it true that for any quotient ${\mathbb S}_{1}$ of ${\mathbb T}_{1}$ there exists a quotient ${\mathbb S}_{2}$ of ${\mathbb T}_{2}$ such that the Morita-equivalence between ${\mathbb T}_{1}$ and ${\mathbb T}_{2}$ restricts to a Morita-equivalence between ${\mathbb S}_{1}$ and ${\mathbb S}_{2}$? The duality theorem gives a straight positive answer to this question. Indeed, the subtoposes of the classifying topos of ${\mathbb T}_{1}$ correspond via the duality theorem, on one hand, to the quotients of ${\mathbb T}_{1}$ and on the other hand, since the classifying topos of ${\mathbb T}_{1}$ is equivalent to the classifying topos of ${\mathbb T}_{2}$, to the quotients of ${\mathbb T}_{2}$. Notice the role of the classifying topos as a bridge (where the invariant is in this case the notion of subtopos and the two different sites of definition of the classifying topos are the geometric syntactic sites of the two theories). The kind of insight brought by the duality theorem in this context is not only theoretical; as we have already remarked, the bijection between quotients and subtoposes can often be exploited for obtaining explicit axiomatizations for geometric theories, and in particular an explicit axiomatization of ${\mathbb S}_{2}$ starting with an axiomatization of ${\mathbb S}_{1}$.
     
Let us now come back to the fundamental role of the duality theorem in allowing the transfer of concepts and results from elementary topos theory into geometric logic.
 
Open (respectively, closed) subtoposes correspond via the duality to quotients obtained by adding sequents of the form $\top \vdash_{[]} \phi$ (respectively, $\phi \vdash_{[]} \bot$), where $\phi$ is a geometric sentence over the signature of the theory, and the surjection-inclusion factorization of a geometric morphism has the following natural semantic interpretation.

\begin{theorem}\label{LatticesTheories_surinclthm}
Let $\mathbb T$ be a geometric theory over a signature $\Sigma$ and $f:{\cal F}\to {\cal E}$ be a geometric morphism into the classifying topos $\cal E$ for $\mathbb T$, corresponding to a $\mathbb T$-model $M$ in $\cal F$ as above. Then the topos ${\cal E}'$ in the surjection-inclusion factorization ${\cal F}\epi {\cal E}' \hookrightarrow {\cal E}$ of $f$ corresponds via the duality theorem to the quotient $Th(M)$ of $\mathbb T$ consisting of all the geometric sequents $\sigma$ over $\Sigma$ which hold in $M$.
\end{theorem} 

A particularly interesting topos-theoretic invariant which, via the duality theorem, leads to important insights on geometric theory is the notion of \emph{Booleanization} of a topos i.e. the subtopos of $\neg\neg$-sheaves. Given a geometric theory $\mathbb T$, we call the quotient of $\mathbb T$ corresponding via the duality theorem to the Booleanization of the classifying topos of $\mathbb T$ the \emph{Booleanization} of $\mathbb T$.

The following result (Theorem 6.4.7 \cite{OC}) gives an explicit axiomatization of the Booleanization of a geometric theory.

\begin{theorem}
Let $\mathbb T$ be a geometric theory over a signature $\Sigma$. Then the Booleanization of $\mathbb T$ is the theory obtained by adding to the axioms of $\mathbb T$ all the geometric sequents of the form $\top \vdash_{\vec{x}} \phi(\vec{x})$, where $\phi(\vec{x})$ is a geometric formula over $\Sigma$ such that for any geometric formula $\chi(\vec{x})$ over $\Sigma$, $\phi \wedge \chi \vdash_{\vec{x}} \bot$ provable in $\mathbb T$ implies $\chi \vdash_{\vec{x}} \bot$ provable in $\mathbb T$. 
\end{theorem}

Booleanizations often specialize, in the case of important mathematical theories, to interesting quotients. For example, the Booleanization of the theory of linear orders is the theory of dense linear orders without endpoints, the Booleanization of the theory of Boolean algebras is the theory of atomless Boolean algebras (cfr. chapter 9 of \cite{OC}), the Booleanization of the (coherent) theory of fields is the (geometric) theory of algebraically closed fields of finite characteristic in which every element is algebraic over the prime field (cfr. \cite{OC4}). Moreover, under appropriate hypotheses, they axiomatize weakly homogeneous models in the sense of classical model theory (cfr. sections \ref{presheaf} and \ref{fraisse} below).  

Another topos-theoretic invariant that has already proved to be relevant for classical Mathematics is the notion of \emph{DeMorganization} of a topos. This invariant was introduced in \cite{OC}, where it was shown that every elementary topos has a largest dense subtopos satisfying De Morgan's law; in the same context, the name \emph{DeMorganization} was used to denote this subtopos. In chapter 6 of \cite{OC}, an explicit axiomatization for the quotient of a geometric theory corresponding via the duality theorem to the DeMorganization of its classifying topos was obtained, and in \cite{OC4}, the authors proved that the DeMorganization of the (coherent) theory of fields is the geometric theory of fields of finite characteristic in which every element is algebraic over the prime field.     

These two examples show that these logically-motivated invariants of toposes have interesting manifestations in a variety of mathematical contexts; the fact that they specialize to important mathematical theories is a clear indication of their centrality in Mathematics.

\newpage   
   
\section{Examples}

In this section we discuss further examples from \cite{OC} which illustrate the application of our philosophy `toposes as bridges'.

\subsection{Theories of presheaf type}\label{presheaf}

A geometric theory is said to be of \emph{presheaf type} to be a geometric theory whose classifying topos is equivalent to a presheaf topos. In other words, a geometric theory is of presheaf type if and only if it is Morita-equivalent to a theory of flat functors on a small category $\cal C$. The class of theories of presheaf type is interesting for many reasons. One important conceptual reason is that any small category $\cal C$ can be regarded, up to Cauchy-completion, as the category $\textrm{f.p.} {\mathbb T}\textrm{-mod}(\Set)$ of (representatives of isomorphism classes of) finitely presentable models of a theory of presheaf type $\mathbb T$ (take $\mathbb T$ to be the theory of flat functors on $\cal C$); as a consequence, any Grothendieck topos is (up to equivalence) of the form $\Sh(\textrm{f.p.} {\mathbb T}\textrm{-mod}(\Set)^{\textrm{op}}, J)$ for some theory of presheaf type $\mathbb T$. Another reason is that this class of theories contains all the cartesian, and in particular all the finitary algebraic, theories, and also many other interesting mathematical theories (for example, the theory of linear orders, cfr. \cite{MM} and \cite{El} and the geometric theory of finite sets, cfr. \cite{El}).  

The classifying topos of a theory of presheaf type $\mathbb T$ is always given by the functor category $[\textrm{f.p.} {\mathbb T}\textrm{-mod}(\Set), \Set]$, where $\textrm{f.p.} {\mathbb T}\textrm{-mod}(\Set)$ is the category of finitely presentable models of $\mathbb T$ (cfr. chapter 4 of \cite{OC}). Therefore, by the duality theorem, the quotients of $\mathbb T$ correspond bijectively to the Grothendieck topologies on the opposite of the category $\textrm{f.p.} {\mathbb T}\textrm{-mod}(\Set)$; we will refer to the Grothendieck topology corresponding to a quotient ${\mathbb T}'$ as to the Grothendieck topology \emph{associated} to ${\mathbb T}'$. This naturally leads to the following question: can the models of a quotient ${\mathbb T}'$ of $\mathbb T$ (in any Grothendieck topos) be characterized among the models of $\mathbb T$ directly in terms of the associated Grothendieck topology? This question was answered affirmatively in Chapter 4 of \cite{OC}; specifically, a notion of \emph{$J$-homogeneous} model of a theory of presheaf type $\mathbb T$ in a Grothendieck topos was introduced (for a Grothendieck topology $J$ on $\textrm{f.p.} {\mathbb T}\textrm{-mod}(\Set)^{\textrm{op}}$), and it was shown that, in every Grothendieck topos, the models of a quotient ${\mathbb T}'$ of $\mathbb T$ are exactly the $J$-homogeneous ones, where $J$ is the Grothendieck topology associated to ${\mathbb T}'$. If the category $\textrm{f.p.} {\mathbb T}\textrm{-mod}(\Set)^{\textrm{op}}$ satisfies the right Ore condition and $J_{at}$ is the atomic topology on it then the notion of $J_{at}$-homogeneous model specializes in $\Set$ to the notion of (weakly) homogeneous model in classical Model Theory (cfr. section \ref{fraisse}).

The investigation of the class of theories of presheaf type provided a splendid opportunity for the author of \cite{OC} to test the methodologies based on her view `toposes as bridges' described above. Indeed, the notion of theory of presheaf type gives automatically rise to a Morita-equivalence; specifically, if $\mathbb T$ is a theory of presheaf type then we have a Morita-equivalence between $\mathbb T$ and the theory of flat functors on the category $\textrm{f.p.} {\mathbb T}\textrm{-mod}(\Set)^{\textrm{op}}$. In topos-theoretic terms, we have an equivalence of Grothendieck toposes $\Sh({\cal C}_{\mathbb T}, J_{\mathbb T}) \simeq [\textrm{f.p.} {\mathbb T}\textrm{-mod}(\Set), \Set]$. The two sites of definitions for this classifying topos are `different enough' for a non-trivial transfer of information from one to another to be attainable by using our methods. In fact, 
in \cite{OC} many results of this kind were obtained, culminating in the topos-theoretic interpretation of Fra\"iss\'e's construction in Model Theory (cfr. section \ref{fraisse}). We shall now briefly survey some of them with the purpose of illustrating how they arise from the application of our principles.

We have already discussed implications of the double representation
\[
\Sh({\cal C}_{\mathbb T}, J_{\mathbb T}) \simeq [\textrm{f.p.} {\mathbb T}\textrm{-mod}(\Set), \Set]
\]
of the classifying topos of $\mathbb T$ for the axiomatization of the quotients of $\mathbb T$. 

The following result provides a link between syntactic properties of a quotient ${\mathbb T}'$ of $\mathbb T$ and `topological' properties of the associated Grothendieck topology $J$ on $\textrm{f.p.} {\mathbb T}\textrm{-mod}(\Set)^{\textrm{op}}$. 
   
\begin{theorem}\label{Onetoposmany_sem}
Let $\mathbb T$ be a theory of presheaf type over a signature $\Sigma$, ${\mathbb T}'$ be a quotient of $\mathbb T$ with associated Grothendieck topology $J$ on $\textrm{f.p.} {\mathbb T}\textrm{-mod}(\Set)^{\textrm{op}}$ and $\phi(\vec{x})$ be a geometric formula over $\Sigma$ which presents a $\mathbb T$-model $M$. Then

\begin{enumerate}[(i)]

\item $\phi(\vec{x})$ is $\mathbb T$-irreducible; in particular, $\phi(\vec{x})$ is $\mathbb T$-provably equivalent to a regular formula;

\item if the site $(\textrm{f.p.} {\mathbb T}\textrm{-mod}(\Set)^{\textrm{op}}, J)$ is locally connected (for example if $\textrm{f.p.} {\mathbb T}\textrm{-mod}(\Set)^{\textrm{op}}$ satisfies the right Ore condition and every $J$-covering sieve is non-empty) then $\phi(\vec{x})$ is ${\mathbb T}'$-indecomposable;
  
\item if $(\textrm{f.p.} {\mathbb T}\textrm{-mod}(\Set)^{\textrm{op}}$ satisfies the right Ore condition and $J$ is the atomic topology on $(\textrm{f.p.} {\mathbb T}\textrm{-mod}(\Set)^{\textrm{op}}$ then $\phi(\vec{x})$ is ${\mathbb T}'$-complete;

\item if every $J$-covering sieve on $M$ contains a $J$-covering sieve generated by a finite family of morphisms then $\phi(\vec{x})$ is ${\mathbb T}'$-compact.

\end{enumerate}
\end{theorem}

The proof of these results consists in the transfer of (properties of) a particular invariant, namely the interpretation of the formula $\phi(\vec{x})$ in the universal model of ${\mathbb T}'$, across the two different representations of the classifying topos of ${\mathbb T}'$, as the category $\Sh(\textrm{f.p.} {\mathbb T}\textrm{-mod}(\Set)^{\textrm{op}}, J)$ and as the category of sheaves $\Sh({\cal C}_{{\mathbb T}'}, J_{{\mathbb T}'})$ on the geometric syntactic site of ${\mathbb T}'$. Specifically, if
\[ 
l^{\textrm{f.p.} {\mathbb T}\textrm{-mod}(\Set)^{\textrm{op}}}_{J}:\textrm{f.p.} {\mathbb T}\textrm{-mod}(\Set)\to \Sh(\textrm{f.p.} {\mathbb T}\textrm{-mod}(\Set)^{\textrm{op}}, J)
\]
is the composite of the associated sheaf functor 
\[
[\textrm{f.p.} {\mathbb T}\textrm{-mod}(\Set), \Set] \rightarrow \Sh(\textrm{f.p.} {\mathbb T}\textrm{-mod}(\Set)^{\textrm{op}}, J)
\]
with the Yoneda embedding 
\[
\textrm{f.p.} {\mathbb T}\textrm{-mod}(\Set)^{\textrm{op}} \rightarrow [\textrm{f.p.} {\mathbb T}\textrm{-mod}(\Set), \Set]
\]
and 
\[
y:{\cal C}_{{\mathbb T}'}\hookrightarrow \Sh({\cal C}_{{\mathbb T}'}, J_{{\mathbb T}'})
\]
is the Yoneda embedding then there is an equivalence 
\[
\Sh({\cal C}_{{\mathbb T}'}, J_{{\mathbb T}'}) \simeq \Sh(\textrm{f.p.} {\mathbb T}\textrm{-mod}(\Set)^{\textrm{op}}, J)
\]
which sends $y(\{\vec{x}. \phi\})$ to $l^{\textrm{f.p.} {\mathbb T}\textrm{-mod}(\Set)^{\textrm{op}}}_{J}(M)$, where $M$ is the $\mathbb T$-model presented by $\phi(\vec{x})$. Now, the transfer takes place as follows: properties of the site $(\textrm{f.p.} {\mathbb T}\textrm{-mod}(\Set)^{\textrm{op}}, J)$ imply properties of the object $l^{\textrm{f.p.} {\mathbb T}\textrm{-mod}(\Set)^{\textrm{op}}}_{J}(M)$ in the classifying topos $\Sh(\textrm{f.p.} {\mathbb T}\textrm{-mod}(\Set)^{\textrm{op}}, J)$, which are transferred via the equivalence $\Sh({\cal C}_{{\mathbb T}'}, J_{{\mathbb T}'}) \simeq \Sh(\textrm{f.p.} {\mathbb T}\textrm{-mod}(\Set)^{\textrm{op}}, J)$ to properties of the object $y(\{\vec{x}. \phi\})$ in the classifying topos $\Sh({\cal C}_{{\mathbb T}'}, J_{{\mathbb T}'})$, which are in turn translated into properties of the syntactic site $({\cal C}_{{\mathbb T}'}, J_{{\mathbb T}'})$ and hence of the theory ${\mathbb T}'$.     

The following result, which expresses the fact that the syntactic and semantic notions of finite presentability of a model (cfr. \cite{OC} for a formal definition of these concepts) coincide for theories of presheaf type is obtained by the same method, by using the invariant property of an object of a Grothendieck topos to be irreducible.

\begin{theorem}\label{Onetoposmany_presheafcomplete}
Let $\mathbb T$ be a theory of presheaf type over a signature $\Sigma$. Then

\begin{enumerate}

\item any finitely presentable $\mathbb T$-model in $\Set$ is presented by a $\mathbb T$-irreducible geometric formula $\phi(\vec{x})$ over $\Sigma$;

\item conversely, any $\mathbb T$-irreducible geometric formula $\phi(\vec{x})$ over $\Sigma$ presents a $\mathbb T$-model.
\end{enumerate}
In particular, the category $\textrm{f.p.} {\mathbb T}\textrm{-mod}(\Set)^{\textrm{op}}$ is equivalent to the full subcategory of ${\cal C}_{\mathbb T}^{\textrm{geom}}$ on the $\mathbb T$-irreducible formulae. 
\end{theorem}

Another interesting result about theories of presheaf type is the following definability theorem (Corollary 7.2.2 \cite{OC}).

\begin{theorem}
Let $\mathbb T$ be a theory of presheaf type over a signature $\Sigma$, let $A_{1}, \ldots, A_{n}$ a string of sorts of $\Sigma$ and suppose we are given, for every finitely presentable $\Set$-model $M$ of $\mathbb T$ a subset $R_{M}$ of $MA_{1}\times \cdots \times MA_{n}$ in such a way that each $\mathbb T$-model homomorphism $h:M\to N$ maps $R_{M}$ into $R_{N}$. Then there exists a geometric formula-in-context $\phi(x^{A_{1}}, \ldots, x^{A_{n}})$ such that $R_{M}=[[\phi]]_{M}$ for each $M$.       
\end{theorem}

This result again is proved by the method of transferring an invariant across different sites of definition of a given classifying topos. In this case, the invariant is the notion of subobject of the underlying object(s) of the universal model of $\mathbb T$ and the two representations of the classifying topos of $\mathbb T$ are $[\textrm{f.p.} {\mathbb T}\textrm{-mod}(\Set)^{\textrm{op}}, \Set]$ and $\Sh({\cal C}_{\mathbb T}, J_{\mathbb T})$. The argument goes as follows. In terms of the representation $[\textrm{f.p.} {\mathbb T}\textrm{-mod}(\Set)^{\textrm{op}}, \Set]$, the universal model of $\mathbb T$ can be described as the functor $M_{\mathbb T}$ in $[\textrm{f.p.} {\mathbb T}\textrm{-mod}(\Set), \Set]$ which assigns to a sort $A$ the functor $M_{\mathbb T}A$ given by $(M_{\mathbb T}A)(M)=MA$ and is defined for function and relation symbols in the obvious way, while in terms of the other representation $\Sh({\cal C}_{\mathbb T}, J_{\mathbb T})$ the universal model of $\mathbb T$ admits a syntactic description such that the interpretation of a formula $\phi(\vec{x})$ in it is given by $y(\{\vec{x}. \phi\})$ (cfr. chapter 7 of \cite{OC} for more details). Therefore, the theorem follows at once from the remark that the assignment $M\to R_{M}$ in the statement of the theorem gives rise to a subobject $R\mono M_{\mathbb T}$ in the topos $[\textrm{f.p.} {\mathbb T}\textrm{-mod}(\Set), \Set]$. 

We note that this result is intimately tied to the peculiar nature of theories of presheaf type. We cannot expect the same theorem to hold for larger classes of theories; for example, the theorem does not hold in general for coherent theories. Indeed, the theory of decidable rings (i.e. the theory of commutative rings with unit equipped with a predicate which is provably complemented to the equality relation) clearly satisfies the hypotheses of the theorem for both the property of an element to be nilpotent and for the complement of such property, and if this latter property were definable by a  geometric formula then each of the properties would be definable by a coherent formula, and this can be easily proved (by means of a compactness argument) not to be true.  

We remark that these results which, as we have seen, naturally follow from the implementation of our idea of `toposes as bridges', are by no means trivial and, while they specialize to useful results in a variety of mathematical contexts, they could be very difficult to prove by using traditional techniques. On the other hand, an infinite number of other results can be derived almost `automatically' by applying the same principles to different toposes and different invariants; the amount of non-trivial mathematical results that can be `automatically generated' by the application of these methodologies is simply boundless. These methods represent a new way of doing mathematics which differs, both formally and substantially, from that given by the classical techniques.

\subsection{Fra\"iss\'e's construction from a topos-theoretic perspective}\label{fraisse}

In this section, we discuss a particularly significant example (from Chapter 9 of \cite{OC}) of a result which, although being obtained as an application of our purely topos-theoretic methods, is directly related to a well-known piece of Mathematics, namely Fra\"iss\'e's construction in Model Theory.   

In section \ref{presheaf} we have discussed applications of the philosophy `toposes as bridges' in the context of invariants of \emph{objects} of toposes;  instead, the main result of this section arises from the consideration of `global' invariants on a given classifying topos.

Before we can state the theorem, we need to give some definitions.

\begin{definition}
A category $\cal C$ is said to satisfy the \emph{amalgamation property} if for every objects $a,b,c\in {\cal C}$ and morphisms $f:a\rightarrow b$, $g:a\rightarrow c$ in $\cal C$ there exists an object $d\in \cal C$ and morphisms $f':b\rightarrow d$, $g':c\rightarrow d$ in $\cal C$ such that $f'\circ f=g'\circ g$:
\[  
\xymatrix {
a \ar[d]_{g} \ar[r]^{f} & b  \ar@{-->}[d]^{f'} \\
c \ar@{-->}[r]_{g'} & d } 
\] 
\end{definition} 

Notice that $\cal C$ satisfies the amalgamation property if and only if ${\cal C}^{\textrm{op}}$ satisfies the right Ore condition. So if $\cal C$ satisfies the amalgamation property then we may equip ${\cal C}^{\textrm{op}}$ with the atomic topology. 

\begin{definition}
A category $\cal C$ is said to satisfy the \emph{joint embedding property} if for every pair of objects $a,b\in {\cal C}$ there exists an object $c\in \cal C$ and morphisms $f:a\rightarrow c$, $g:b\rightarrow c$ in $\cal C$:
\[  
\xymatrix {
 & a  \ar@{-->}[d]^{f} \\
b \ar@{-->}[r]_{g} & c } 
\] 
\end{definition} 

\begin{definition}
Let $\mathbb T$ be a theory of presheaf type. A model $M$ of $\mathbb T$ is said to be \emph{homogeneous} if for any models $a,b \in \textsf{f.p.} {\mathbb T}\textsf{-mod}(\Set)$ and arrows $j:a\rightarrow b$ and $\chi:a\rightarrow M$ in ${\mathbb T}\textsf{-mod}(\Set)$ there exists an arrow $\tilde{\chi}:b\rightarrow M$ in ${\mathbb T}\textsf{-mod}(\Set)$ such that $\tilde{\chi}\circ j=\chi$:
\[  
\xymatrix {
a \ar[d]_{j} \ar[r]^{\chi} & M \\
b \ar@{-->}[ur]_{\tilde{\chi}} &  } 
\] 
\end{definition}
  
Our main result is the following. 

\begin{theorem}\label{AtomicTheories_teoFra\"iss\'e}
Let $\mathbb{T}$ be a theory of presheaf type such that the category $\textrm{f.p.} {\mathbb T}\textrm{-mod}(\Set)$ satisfies the amalgamation and joint embedding properties. Then any two countable homogeneous $\mathbb T$-models in $\Set$ are isomorphic.
\end{theorem}

We can describe the structure of the proof of this theorem as follows.

If the category $\textrm{f.p.} {\mathbb T}\textrm{-mod}(\Set)$ is empty then $\mathbb T$ has no models in $\Set$ and hence the thesis is trivially true. We will therefore assume that $\textrm{f.p.} {\mathbb T}\textrm{-mod}(\Set)$ is non-empty. Since $\textrm{f.p.} {\mathbb T}\textrm{-mod}(\Set)$ satisfies the amalgamation property, we can equip the opposite category $\textrm{f.p.} {\mathbb T}\textrm{-mod}(\Set)^{\textrm{op}}$ with the atomic topology $J_{at}$. Our argument is based on the topos $\Sh(\textrm{f.p.} {\mathbb T}\textrm{-mod}(\Set)^{\textrm{op}}, J_{at})$. By the duality theorem, the subtopos 
\[
\Sh(\textrm{f.p.} {\mathbb T}\textrm{-mod}(\Set)^{\textrm{op}}, J_{at}) \hookrightarrow [\textrm{f.p.} {\mathbb T}\textrm{-mod}(\Set)^{\textrm{op}}, \Set]
\]
corresponds to a unique quotient (up to syntactic equivalence) ${\mathbb T}'$ of $\mathbb T$, such that its classifying topos is $\Sh(\textrm{f.p.} {\mathbb T}\textrm{-mod}(\Set)^{\textrm{op}}, J_{at})$ (Note that, conceptually, the use of the duality theorem essentially amounts to transferring the invariant notion of subtopos from the representation $[\textrm{f.p.} {\mathbb T}\textrm{-mod}(\Set)^{\textrm{op}}, \Set]$ of the classifying topos of $\mathbb T$ to the syntactic representation $\Sh({\cal C}_{\mathbb T}, J_{\mathbb T})$ of it). From the technique of `Yoneda representation of flat functors' (cfr. Chapter 4 of \cite{OC}), it follows that in every Grothendieck topos, in particular in $\Set$, the models of ${\mathbb T}'$ are exactly the homogeneous $\mathbb T$-models. 

Now, the `core' of the argument consists in transferring the invariant property of a topos to be \emph{atomic and two-valued} topos from the `semantic' representation $\Sh(\textrm{f.p.} {\mathbb T}\textrm{-mod}(\Set)^{\textrm{op}}, J_{at})$ to the syntactic representation $\Sh({\cal C}_{{\mathbb T}'}, J_{{\mathbb T}'})$ of the classifying topos of ${\mathbb T}'$. One can prove that if $\textrm{f.p.} {\mathbb T}\textrm{-mod}(\Set)$ satisfies the joint embedding property, as well as the amalgamation property, then the topos $\Sh(\textrm{f.p.} {\mathbb T}\textrm{-mod}(\Set)^{\textrm{op}}, J_{at})$ is two-valued, as well as atomic (in fact, the other implication also holds under our hypotheses). But if we read this invariant from the point of view of the other representation $\Sh({\cal C}_{{\mathbb T}'}, J_{{\mathbb T}'})$, this tells us precisely that the theory ${\mathbb T}'$ is `atomic' and `complete', which forces ${\mathbb T}'$ to be countably categorical (cfr. Chapter 8 \cite{OC}). In view of the aforementioned fact that the models of ${\mathbb T}'$ are exactly the homogeneous $\mathbb T$-models, this concludes the proof of the theorem. By transferring the invariant construction of  the Booleanization across the two different representations of the classifying topos of $\mathbb T$, one immediately sees that the quotient ${\mathbb T}'$ is precisely the Booleanization of $\mathbb T$.
  
Another interesting application of our methodology of transferring invariants across different sites of definition of one topos (given in Chapter 9 of \cite{OC}) concerns the existence of homogeneous models of a theory of presheaf type $\mathbb T$ in $\Set$. As we have already remarked, there is a quotient ${\mathbb T}'$ of $\mathbb T$ which axiomatizes the homogeneous models in every Grothendieck topos and is classified by the topos $\Sh(\textrm{f.p.} {\mathbb T}\textrm{-mod}(\Set)^{\textrm{op}}, J_{at})$. Therefore, the existence of a homogeneous model of $\mathbb T$ in $\Set$ is equivalent to the existence of a point of the topos $\Sh(\textrm{f.p.} {\mathbb T}\textrm{-mod}(\Set)^{\textrm{op}}, J_{at})$. Notice that this latter property is a topos-theoretic invariant. Now, one can easily prove that if a topos $\cal E$ has enough points then $\cal E$ has a point if and only if it is non-trivial; and one can directly prove that the topos $\Sh(\textrm{f.p.} {\mathbb T}\textrm{-mod}(\Set)^{\textrm{op}}, J_{at})$ is non-trivial provided that the category $\textrm{f.p.} {\mathbb T}\textrm{-mod}(\Set)$ is non-empty. Therefore, under the hypotheses that $\textrm{f.p.} {\mathbb T}\textrm{-mod}(\Set)$ should be non-empty, we are reduced to prove that $\Sh(\textrm{f.p.} {\mathbb T}\textrm{-mod}(\Set)^{\textrm{op}}, J_{at})$ has enough points. The interesting thing is that this can be done in many different ways; for example, by Deligne's theorem, it suffices to find a coherent theory which is classified by the topos $\Sh(\textrm{f.p.} {\mathbb T}\textrm{-mod}(\Set)^{\textrm{op}}, J_{at})$. In fact, in \cite{OC}, we exploited the fact that a purely combinatorial condition on the category $\textrm{f.p.} {\mathbb T}\textrm{-mod}(\Set)$ guarantees that the theory of $J_{at}$-continuous flat functors on $\textrm{f.p.} {\mathbb T}\textrm{-mod}(\Set)^{\textrm{op}}$ is coherent.

\subsection{Other examples}\label{examp}

In this section, we discuss other applications of the philosophy `toposes as bridges' to questions in Logic.

Recall that for any fragment of geometric logic, there is a corresponding notion of provability for theories in that fragment; for example, we have a notion of provability of regular sequents in regular logic and a notion of provability of coherent sequents in coherent logic, as well as the classical notion of provability of geometric sequents in geometric logic.

A natural question is whether these notions of provability are compatible with each other, that is if the notion of provability in a given fragment of logic specializes to the notions of provability in a smaller fragment (note that any theory in a given fragment can be seen as a theory in a larger fragment). We show that there is a natural topos-theoretic way of dealing with these problems, which exploits the fact that for theories in a proper fragment of geometric logic one has multiple syntactic representations of their classifying topos, each corresponding to a particular fragment of logic in which the theory lies (cfr. section \ref{onetopos} above). For example, if a theory $\mathbb T$ over a signature $\Sigma$ is coherent then we have two syntactic representations of its classifying topos, as the category of sheaves $\Sh({\cal C}_{\mathbb T}^{\textrm{coh}}, J_{{\cal C}_{\mathbb T}^{\textrm{coh}}})$ on the coherent syntactic site $({\cal C}_{\mathbb T}^{\textrm{coh}}, J_{{\cal C}_{\mathbb T}^{\textrm{coh}}})$ of $\mathbb T$ and as the category of sheaves $\Sh({\cal C}_{\mathbb T}, J_{{\cal C}_{\mathbb T}})$ on the geometric syntactic site $({\cal C}_{\mathbb T}, J_{{\cal C}_{\mathbb T}})$ of $\mathbb T$. From the fact that the two Yoneda embeddings $y_{\textrm{coh}}:{\cal C}_{\mathbb T}^{\textrm{coh}} \to \Sh({\cal C}_{\mathbb T}^{\textrm{coh}}, J_{{\cal C}_{\mathbb T}^{\textrm{coh}}})$ and $y:{\cal C}_{\mathbb T}\to \Sh({\cal C}_{\mathbb T}, J_{{\cal C}_{\mathbb T}})$ are conservative, it follows that a coherent sequent over $\Sigma$ is provable in the universal model of $\mathbb T$ lying in the topos $\Sh({\cal C}_{\mathbb T}^{\textrm{coh}}, J_{{\cal C}_{\mathbb T}^{\textrm{coh}}})$ if and only if it is provable in $\mathbb T$ by using coherent logic, and that a geometric sequent over $\Sigma$ is provable in the universal model of $\mathbb T$ lying in the topos $\Sh({\cal C}_{\mathbb T}^{\textrm{coh}}, J_{{\cal C}_{\mathbb T}^{\textrm{coh}}})$ if and only if it is provable in $\mathbb T$ by using geometric logic. But the property of validity of a given sequent in the universal model of a geometric theory is clearly a topos-theoretic invariant, from which it follows that a coherent sequent is provable in $\mathbb T$ by using coherent logic if and only if it is provable in $\mathbb T$ by using geometric logic, as required.   

Another interesting application concerns `classical completeness' of theories with respect to geometric logic. If $\mathbb T$ is a coherent theory over a signature $\Sigma$, it is well-known that if a coherent sequent over $\Sigma$ is satisfied in every $\Set$-based model of $\mathbb T$ then it is provable in $\mathbb T$ by using coherent logic. It is thus natural to wonder whether, under the same assumption that $\mathbb T$ be coherent, the same holds for geometric sequents i.e. if any geometric sequent over $\Sigma$ which is valid in every $\Set$-based model of $\mathbb T$ is provable in $\mathbb T$ by using geometric logic. The answer to this question is affirmative, and naturally follows from the transfer of an invariant, namely the property of a topos to have enough points, across the two different representations $\Sh({\cal C}_{\mathbb T}^{\textrm{coh}}, J_{{\cal C}_{\mathbb T}^{\textrm{coh}}}) \simeq \Sh({\cal C}_{\mathbb T}, J_{{\cal C}_{\mathbb T}})$ of the classifying topos of $\mathbb T$. Indeed, in terms of the site $({\cal C}_{\mathbb T}^{\textrm{coh}}, J_{{\cal C}_{\mathbb T}^{\textrm{coh}}})$ the invariant rephrases as the condition that every coherent sequent over $\Sigma$ which is satisfied in all the $\Set$-based models of $\mathbb T$ should be provable in $\mathbb T$ by using coherent logic, while in terms of the site $({\cal C}_{\mathbb T}, J_{{\cal C}_{\mathbb T}})$ it rephrases as the condition that every geometric sequent over $\Sigma$ which is satisfied in all the $\Set$-based models of $\mathbb T$ should be provable in $\mathbb T$ by using geometric logic.

Another result which follows quite naturally from the double representation $\Sh({\cal C}_{\mathbb T}^{\textrm{coh}}, J_{{\cal C}_{\mathbb T}^{\textrm{coh}}}) \simeq \Sh({\cal C}_{\mathbb T}, J_{{\cal C}_{\mathbb T}})$ of the classifying topos of a coherent theory is the following (Theorem 10.2.5 \cite{OC}).

\begin{theorem}\label{Onetoposmany_cohcrit}
Let $\mathbb T$ be a geometric theory over a signature $\Sigma$. Then $\mathbb T$ is coherent if and only if for any coherent formula $\{\vec{x}. \phi\}$ over $\Sigma$, for any family $\{\psi_{i}(\vec{x}) \textrm{ | } i\in I\}$ of geometric formulae in the same context, $\phi \vdash_{\vec{x}} \mathbin{\mathop{\textrm{\huge $\vee$}}\limits_{i\in I}} \psi_{i}$ provable in $\mathbb T$ implies $\phi \vdash_{\vec{x}} \mathbin{\mathop{\textrm{\huge $\vee$}}\limits_{i\in I'}} \psi_{i}$ provable in $\mathbb T$ for some finite subset $I'$ of $I$. 
\end{theorem}

The non-trivial part of the theorem follows immediately from the remark that if $\mathbb T$ is a coherent theory over a signature $\Sigma$ we have an equivalence of classifying toposes $\Sh({\cal C}_{\mathbb T}^{\textrm{coh}}, J_{{\cal C}_{\mathbb T}^{\textrm{coh}}}) \simeq \Sh({\cal C}_{\mathbb T}, J_{{\cal C}_{\mathbb T}})$ such that for any coherent formula $\{\vec{x}. \phi\}$ over $\Sigma$, the objects $y_{\textrm{coh}}(\{\vec{x}.\phi\})$ and  $y(\{\vec{x}.\phi\})$ correspond to each other under the equivalence (where $y_{\textrm{coh}}$ and $y$ are the two Yoneda embeddings $y_{\textrm{coh}}:{\cal C}_{\mathbb T}^{\textrm{coh}} \to \Sh({\cal C}_{\mathbb T}^{\textrm{coh}}, J_{{\cal C}_{\mathbb T}^{\textrm{coh}}})$ and $y:{\cal C}_{\mathbb T}\to \Sh({\cal C}_{\mathbb T}, J_{{\cal C}_{\mathbb T}})$). If we consider the invariant property of an object of a topos to be compact, we see that the fact that the site $({\cal C}_{\mathbb T}^{\textrm{coh}}, J_{{\cal C}_{\mathbb T}^{\textrm{coh}}})$ is coherent implies that the object $y_{\textrm{coh}}(\{\vec{x}.\phi\})$ is compact, but from the point of view of the geometric syntactic site $({\cal C}_{\mathbb T}, J_{{\cal C}_{\mathbb T}})$ the property of $y(\{\vec{x}.\phi\})$ to be compact reprases as follows: for any family $\{\psi_{i}(\vec{x}) \textrm{ | } i\in I\}$ of geometric formulae in the same context, $\phi \vdash_{\vec{x}} \mathbin{\mathop{\textrm{\huge $\vee$}}\limits_{i\in I}} \psi_{i}$ provable in $\mathbb T$ implies $\phi \vdash_{\vec{x}} \mathbin{\mathop{\textrm{\huge $\vee$}}\limits_{i\in I'}} \psi_{i}$ provable in $\mathbb T$ for some finite subset $I'$ of $I$.  
  
This theorem shows that what characterizes coherent theories among geometric ones is precisely a general form of compactness. 

Of course, all of the above-mentioned results have a version for cartesian or regular theories in place of coherent theories (cfr. chapter 10 of \cite{OC}).

Finally, we note that the well-known equivalence between Deligne's theorem in Algebraic Geometry asserting that every coherent topos has enough points and the classical completeness theorem for coherent theories in Logic can be naturally expressed in our framework. Indeed, starting from a coherent topos represented as the category of sheaves on a site $({\cal C}, J)$ such that $\cal C$ has finite limits and $J$ is a finite-type Grothendieck topology on it, one can consider a different site of definition of it, namely the coherent syntactic site $({\cal C}_{\mathbb T}^{\textrm{coh}}, J_{{\cal C}_{\mathbb T}^{\textrm{coh}}})$ of a coherent theory $\mathbb T$ (since, under our hypotheses, the theory of flat $J$-continuous functors on $\cal C$ admits a coherent axiomatization). So one has two different representations for the given coherent topos, one of `geometric' nature and another one of logical nature; and the above-mentioned equivalence follows precisely from the fact that the expression of the invariant property of having enough points in terms of the coherent syntactic site $({\cal C}_{\mathbb T}^{\textrm{coh}}, J_{{\cal C}_{\mathbb T}^{\textrm{coh}}})$ of $\mathbb T$ amounts precisely to the classical completeness theorem for $\mathbb T$.

\section{Toposes for the working mathematician}\label{work}

In this section, we summarize some of the aspects of Topos Theory, and in particular of the methodologies that we have described above in this paper, which should make the subject attractive to the `working mathematician'. 

We have seen that to every first-order mathematical theory (admitting a geometric axiomatization) one can naturally associate a classifying topos, which extracts the essential features of the theory, that is those features which are invariant with respect to Morita-equivalence. A given topos can have many different sites of definition, which illuminate different aspects of the theories classified by it. Topos-theoretic invariants (i.e. properties of toposes or constructions involving them which are invariant under categorical equivalence) can then be used to transfer properties between theories classified by that topos. In fact, the abstract relationship between a site $({\cal C},J)$ and the topos $\Sh({\cal C}, J)$ which it `generates' is often very natural, enabling us to easily transfer invariants across different sites. Moreover, as we have seen in sections \ref{logic} and \ref{duality}, the level of generality represented by topos-theoretic invariants is ideal to capture several important features of mathematical theories.        

The kind of `coding' $({\cal C}, J)\to \Sh({\cal C}, J)$ given by the operation of taking the category of sheaves on a given site is highly non-trivial; the two ingredients $({\cal C}, J)$ are combined together in a sophisticated way to form the topos $\Sh({\cal C}, J)$. The role of topos-theoretic invariants is precisely that of `decoding' the information hidden in a topos, in the sense of identifying properties of a site $({\cal C}, J)$ which are implied by properties of the corresponding topos $\Sh({\cal C}, J)$, or describing the result of constructions defined on toposes directly in terms of the site in a way which is intelligible to specialists outside Topos Theory. Indeed, while the notion of site is quite near to the common mathematical practice, a topos is really an object of a different nature, a peculiar entity which can be studied by means of specific methods (i.e. those of Topos Theory) which are \emph{not} those of `classical Mathematics' (cfr. section \ref{genetics} for a conceptual analogy which may illuminate this point). 

The link between classical mathematics and Topos Theory is given by the fact that many constructions or properties naturally arising in Mathematics can be expressed in terms of topos-theoretic invariants; for example, cohomology groups of toposes are topos-theoretic invariants that specialize, in the case of particular toposes, to cohomologies arising in mathematical practice (such as for example the usual cohomology of topological spaces, Galois cohomology, Eilenberg-Mac Lane cohomology, Weil cohomologies), allowing at the same time enough freedom to construct new `concrete' cohomologies which can be useful for other purposes. 

Now, the advantage of working with toposes rather with other kinds of entities lies in the fact that - as we have seen in the paper - toposes can naturally act as unifying spaces allowing the transfer of properties or constructions between distinct mathematical theories, thus opening the way to the use of methods in one branch of Mathematics to solve problems in another. For example, a property of a cohomology group of a given topos could well be translated, by using a representation of the topos as the category of sheaves on the syntactic site of a geometric theory $\mathbb T$, into some logical property of the theory $\mathbb T$ (similarly, using a site of algebraic nature, one could end up with a property tractable by the methods of algebra, and so on).

The `working mathematician' could very well attempt to formulate his or her properties of interest in terms of topos-theoretic invariants, and derive equivalent versions of them by using alternative sites. Also, he or she could try to express a mathematical duality or equivalence that he or she wants to investigate in terms of a Morita-equivalence, calculate the classifying topos of the two theories and apply topos-theoretic methods to extract new information about it (cfr. sections \ref{onetopos}, \ref{morita} and \ref{bridge}). 

Therefore, much research effort should be done to introduce new invariants for toposes which connect interestingly with questions studied by mathematicians, and to identify Morita-equivalences (also in the form of establishing general representation theorems for toposes); we hope to have convinced the reader that, by only using the most well-known invariants of toposes, one can easily obtain substantial results in distinct mathematical contexts. Site characterizations should then be sought for these invariants in order to allow the transfer of properties to be feasible in presence of Morita-equivalences. 

We remark that a mathematical theory alone can `generate' a considerable number of non-trivial Morita-equivalences (cfr. section \ref{onetopos}); so topos-theoretic methods can be employed to extract information on a single theory, without one having to specifically look for Morita-equivalences connecting that theory to very `different-looking' ones.    
   
In conclusion, Topos Theory has indeed very much to offer to specialists in any field of Mathematics; by means of the methodologies described above in the paper, one can generate a huge number of new results in any mathematical field without any creative effort. This unifying machinery has really the potential to automatically generate results; of course, many of the insights obtained in this way will appear `weird' to the working mathematician (although they might still be quite deep) and in fact, in order to get interesting results, one has to carefully select invariants and Morita-equivalences in such a way that they naturally relate to the questions of interest (as for example in the case of Fra\"iss\'e's construction). Still, one should bear in mind that the kind of insights that these methods can bring into Mathematics is \emph{intrinsically} different from anything that one could naturally expect. It is precisely where traditional methods fail that topos-theoretic methods can possibly lead to a solution (think for example of the proof of the Weil's conjectures by Grothendieck and Deligne).

\subsection{A comparison with Genetics}\label{genetics}

The kind of insight that topos-theoretic methods can bring into Mathematics is comparable to that which Genetics is bringing into Medicine (apologies are due to biologists for any inaccuracies about their subject that may show through the next paragraph). 

As the (human) DNA embodies the essential features of an individual, so the classifying topos embodies the essential features of a mathematical theory. As the DNA can be extracted in many different ways (for example, from different parts of an individual), so the classifying topos can represented and calculated in alternative ways (emphasizing distinct readings of the core of a theory). As the DNA is invariant with respect to the particular physical appearance of the individual at a given time (for example, with respect to age), so the classifying topos is invariant with respect to particular presentations of the theory (for example, with respect to a particular axiomatization of it over its signature). As the DNA of an individual can be studied by using appropriate techniques, which are not those of traditional Medicine, so classifying toposes can be studied by using peculiar methods (i.e. those of Topos Theory) which, although being entirely mathematical, differ from those of classical Mathematics. As in Genetics one studies how modifications of the DNA influence the characteristics of an individual, similarly in Topos Theory one can study the effect that topos-theoretic operations on toposes have on the theories classified by them. As the role of the DNA is that of a unifying concept enabling us to compare individuals with each other, point out differences and discover similarities, so the notion of classifying topos is a unifying one, enabling us to compare distinct mathematical theories with each other and transfer knowledge between them.

We hope to have conveyed by this metaphor the intuition behind the presented results as well as the potential role of topos-theoretic techniques in the Mathematics of the future.

\vspace{0.5cm}
{\flushleft 
{\bf Acknowledgements:}} I am grateful to Marco Benini for his useful comments on a draft of this paper.

\end{document}